\documentclass[12pt,twoside]{article}
\usepackage{amsmath,amsfonts,amssymb,a4,enumerate,epsfig,array,theorem,psfrag}

\newtheorem{theo}{Theorem}

\newtheorem{prop}{Proposition}[section]

\newtheorem{defin}[prop]{Definition}
\newtheorem{lemma}[prop]{Lemma}

\newcommand{\base}{\bullet}

\newcommand{\xx}{{\bf x}}
\newcommand{\ff}{{\bf f}}
\newcommand{\bfg}{{\bf g}}
\newcommand{\bfn}{{\bf n}}

\newcommand{\supp}{\operatorname{supp}}

\newcommand{\interior}{\operatorname{int}}

\newcommand{\ZZ}{{\mathbb{Z}}}

\newcommand{\RR}{{\mathbb{R}}}
\newcommand{\CC}{{\mathbb{C}}}
\newcommand{\Ss}{{\mathbb{S}}}

\newcommand{\DD}{{\mathbb{D}}}

\newcommand{\cF}{{\cal F}}
\newcommand{\cC}{{\cal C}}
\newcommand{\cE}{{\cal E}}
\newcommand{\cEp}{{\cal E}^+}
\newcommand{\cL}{{\cal L}}

\newcommand{\cI}{{\cal I}}
\newcommand{\cJ}{{\cal J}}
\newcommand{\cK}{{\cal K}}
\newcommand{\cO}{{\cal O}}
\newcommand{\cUp}{{\cal U}^+}
\newcommand{\cUpu}{{\cal U}^1}
\newcommand{\cW}{{\cal W}}
\newcommand{\cT}{{\cal T}}
\newcommand{\fF}{{\mathfrak F}}

\newcommand{\nobf}{\noindent\bf}

\def\qed{\unskip\nobreak\hfil\penalty50\hskip1.75em\null\nobreak\hfil
$\blacksquare$ {\parfillskip=0pt \finalhyphendemerits=0 \par}\goodbreak}

\begin{document}
\title{The homotopy and cohomology of spaces \\
of locally convex curves in the sphere --- II}
\author{Nicolau C. Saldanha}
\maketitle

\begin{abstract}
A smooth curve $\gamma: [0,1] \to \Ss^2$ is
locally convex if its geodesic curvature is positive at every point.
J.~A.~Little showed that the space
of all locally positive curves $\gamma$
with $\gamma(0) = \gamma(1) = e_1$ and $\gamma'(0) = \gamma'(1) = e_2$
has three connected components $\cL_{-1,c}$, $\cL_{+1}$, $\cL_{-1,n}$.
The space $\cL_{-1,c}$ is known to be contractible
but the topology of the other two connected components is not well understood.
We prove that all connected components of $\cL_I$ are simply connected,
that $H^2(\cL_{+1};\ZZ) = \ZZ^2$ and $H^2(\cL_{-1,n};\ZZ) = \ZZ$.
\end{abstract}

\section{Introduction}

\footnote{2000 {\em Mathematics Subject Classification}.
Primary 57N65, 53C42; Secondary 34B05.
{\em Keywords and phrases} Convex curves,
topology in infinite dimension, 
periodic solutions of linear ODEs.}

\begin{figure}[ht]
\begin{center}
\epsfig{height=30mm,file=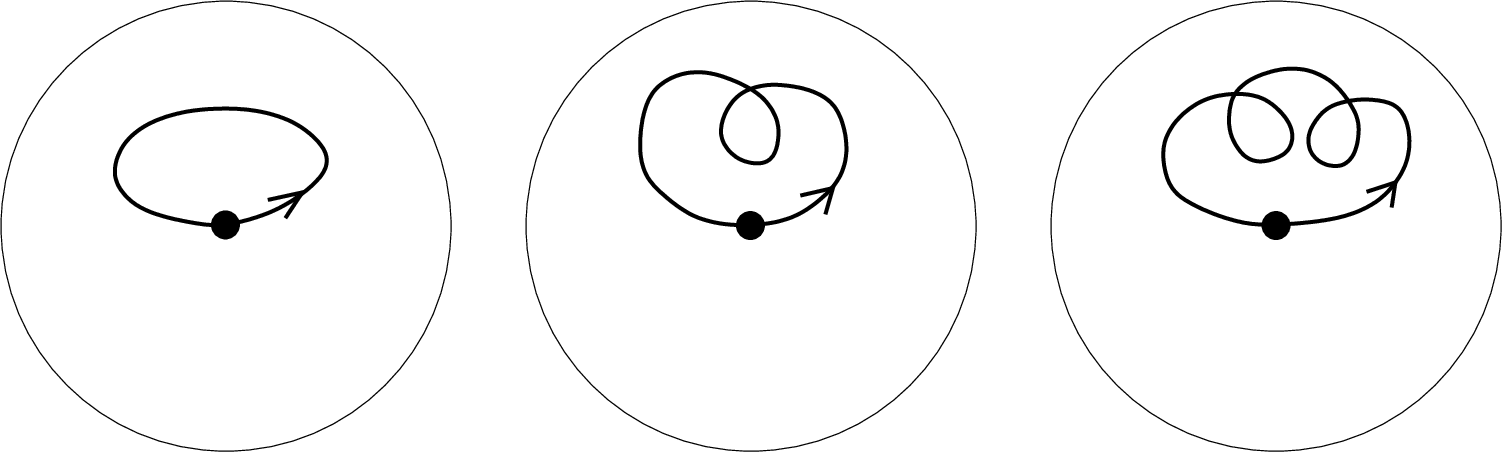}
\end{center}
\caption{Curves in $\cL_{-1,c}$, $\cL_{+1}$ and $\cL_{-1,n}$.}
\label{fig:3comp}
\end{figure}

A curve $\gamma: [0,1] \to \Ss^2$ is called {\sl locally convex}
if its geodesic curvature is always positive,
or, equivalently, if $\det(\gamma(t), \gamma'(t), \gamma''(t)) > 0$ for all $t$.
Let $\cL_I$ be the space of all locally convex curves $\gamma$
with $\gamma(0) = \gamma(1) = e_1$ and $\gamma'(0) = \gamma'(1) = e_2$.
J.~A.~Little \cite{Little} showed that $\cL_I$
has three connected components $\cL_{-1,c}$, $\cL_{+1}$, $\cL_{-1,n}$:
we call these the Little spaces.
Figure \ref{fig:3comp} shows examples of curves in
$\cL_{-1,c}$, $\cL_{+1}$ and $\cL_{-1,n}$, respectively.
The space $\cL_{-1,c}$ is known to be contractible (\cite{ShapiroM})
but the topology of the other two connected components is not well understood.
In this series of papers we present new results concerning
the homotopy and cohomology of the Little spaces.
A more ambitious aim would be to determine the homotopy type
of these spaces (which we hope to accomplish in \cite{Saldanha3}).

Let $\cI_I$ be the space of immersed curves $\gamma: [0,1] \to \Ss^2$,
$\gamma(0) = \gamma(1) = e_1$, $\gamma'(0) = \gamma'(1) = e_2$.
For each $\gamma \in \cI_I$, consider its Frenet frame
$\fF_\gamma: [0,1] \to SO(3)$ and its lift
$\tilde\fF_\gamma: [0, 1] \to \Ss^3$.
The value of $\tilde\fF_\gamma(1)$ defines the two connected components
of $\cI_I$:
$\gamma \in \cI_{+1}$ if and only if $\tilde\fF_\gamma(1) = 1$.
It is well know that each space $\cI_{\pm 1}$ is homotopically equivalent
to $\Omega \Ss^3$.
In particular each $\cI_{\pm 1}$ is connected, simply connected
and there is an element $\xx \in H^2(\cI_{\pm 1};\ZZ)$
such that each $H^{2k}(\cI_{\pm 1};\ZZ) = \ZZ$ is generated by $\xx^k$.
Let $\cL_{+1} = \cI_{+1} \cap \cL_I$,
$\cL_{-1} = \cL_{-1,c} \sqcup \cL_{-1,n} = \cI_{-1} \cap \cL_I$.
In the first paper (\cite{S1}) we saw that the inclusions
$\cL_{\pm 1} \subset \cI_{\pm 1}$ are homotopically surjective
but not homotopy equivalences.
Indeed, we constructed elements $\ff_{2k} \in H^{2k}(\cL_{(-1)^{(k+1)}};\ZZ)$
and maps $\bfg_{2k}: \Ss^{2k} \to \cL_{(-1)^{(k+1)}}$ with
$\ff_{2k}(\bfg_{2k}) = 1$,
$\bfg_{2k}$ homotopic to a constant in $\cI_{(-1)^{(k+1)}}$
and $\ff_{2k}$ not in (the image of) $H^{2k}(\cI_{(-1)^{(k+1)}};\ZZ)$.
In other words, there we give lower estimates for the groups
$H^{2k}(\cL_{\pm 1};\ZZ)$ and $\pi_{2k}(\cL_{\pm 1})$.
In the present paper we give upper estimates which
imply the following theorem.

\begin{theo}
\label{theo:main}
The connected components of $\cL_I$ are simply connected.
Furthermore, $H^2(\cL_{-1,n};\ZZ)$ is generated by $\xx$
and $H^2(\cL_{+1};\ZZ)$ is generated by $\xx$ and $\ff_2$.
Also, $\pi_2(\cL_{-1,n}) = \ZZ$ and $\pi_2(\cL_{+1}) = \ZZ^2$
is generated by $\bfg_2$ and $\tilde\bfg_2$.
\end{theo}

Notice that $\cL_{per}$, the set of all $1$-periodic locally convex
curves $\tilde\gamma: \RR \to \Ss^2$ is homeomorphic to $SO(3) \times \cL_I$:
define $\Psi: \cL_{per} \to SO(3) \times \cL_I$ by
$\Psi(\tilde\gamma) =
(\fF_{\tilde\gamma}(0), (\fF_{\tilde\gamma}(0))^{-1} \tilde\gamma|_{[0,1]})$.
We usually prefer to work in $\cL_I$ but sometimes move to $\cL_{per}$.

In Section 2 we review some known results.
Section 3 contains an algebraic description
of the all-important construction $\Delta^\sharp$;
pulling one loop around is a special case of $\Delta^\sharp$.
Section 4 discusses the uses and limitations of $\Delta^\sharp$
to prove that a map $f: K \to \cL_I$ is homotopic to $\nu_2 \ast f$
(which essentially reduces the problem to the well-understood
scenario of immersions).
In Section 5 the discussion becomes more geometric
and less algebraic as we discuss loops and
the set $\cT_0 \subset \cL_{+1}$ of stars:
roughly, once $\cT_0$ is removed,
$\cL_{+1} \smallsetminus \cT_0$ rather resembles $\cI_{+1}$.
Section 6 polishes a few nasty configurations so that
we can complete the proof of our main results in Section 7.
Section 8 is a very short conclusion.

The author would like to thank Dan Burghelea and
Boris Shapiro for helpful conversations.
The author acknowledges the hospitality of
The Mathematics Department of The Ohio State University
during the winter quarters of 2004 and 2009
and the support of CNPq, Capes and Faperj (Brazil).

\section{Previous results}

\begin{figure}[ht]
\begin{center}
\psfrag{(a)}{(a)}
\psfrag{(b)}{(b)}
\psfrag{(c)}{(c)}
\psfrag{gamma}{$\gamma$}
\psfrag{gammasharp}{$\gamma^\sharp$}
\psfrag{t0}{$t_0$}
\psfrag{t1}{$t_1$}
\epsfig{height=25mm,file=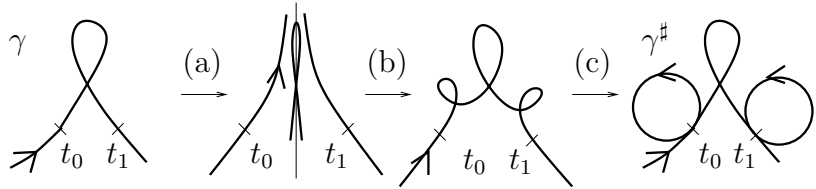}
\end{center}
\caption{How to go from $\gamma$ to $\gamma^\sharp$}
\label{fig:little}
\end{figure}

One of the fundamental constructions in Little's argument
is that if the curve $\gamma \in \cL_I$ has a loop,
we can add a pair of loops as in Figure \ref{fig:little}:
in (a), the loop moves one full turn along a geodesic 
and in (b) the large loops are shrunk (we will discuss (c) later).

Let $\cC_0$ be the circle with diameter $e_1e_3$,
parametrized by $\nu_1 \in \cL_I$,
\[ \nu_1(t) = \left( \frac{1 + \cos(2 \pi t)}{2},
\frac{\sqrt{2}}{2} \sin(2 \pi t),
\frac{1 - \cos(2 \pi t)}{2} \right). \]
For positive $n$, let $\nu_n(t) = \nu_1(nt)$
so that $\nu_1 \in \cL_{-1,c}$ and, for $n > 1$,
$\nu_n \in \cL_{(-1)^n}$.

For $\gamma_1 \in \cI_{\sigma_1}$, $\gamma_2 \in \cI_{\sigma_2}$,
$\sigma_i \in \{+1,-1\}$,
let $\gamma_1 \ast \gamma_2 \in \cI_{\sigma_1 \sigma_2}$
be defined by
\[ (\gamma_1 \ast \gamma_2)(t) = \begin{cases}
\gamma_1(2t), & 0 \le t \le 1/2, \\
\gamma_2(2t-1), & 1/2 \le t \le 1. \end{cases} \]
Notice that if $\gamma_1, \gamma_2 \in \cL_I$
then $\gamma_1 \ast \gamma_2 \in \cL_I$.
For $f: K \to \cI_I$, let $\nu_n \ast f: K \to \cI_I$
be defined by $(\nu_n \ast f)(p) = \nu_n \ast (f(p))$.
Intuitively, $\nu_n \ast f$ is obtained from $f$
by adding $n$ loops to $f(p)$ at the point $f(p)(0)$.
We may want to spread out $n$ loops along the curve:
for $\gamma \in \cI_{\sigma}$ and for large $n$,
define $(F_{n}(\gamma))(t) = \fF_{\gamma}(t) \nu_{n}(t)$.
For small $n$, the above function from $[0,1]$ to $\Ss^2$
may not be an immersion.
For sufficiently large $n$, however,
$F_n(\gamma) \in \cL_{(-1)^n \sigma}$.

By the above construction, it is clear that given
$\gamma \in \cL_{+1} \sqcup \cL_{-1,n}$
there exist $H_a, H_b: [0,1] \to \cL_{\pm 1}$,
$H_a(0) = H_b(0) = \gamma$,
$H_a(1) = \nu_2 \ast \gamma$, $H_b(1) = F_{2n}(\gamma)$.
The construction is not uniform, however:
is depends on the choice of the loop.
In other words, given a compact set $K$
and a map $f: K \to \cL_{+1} \sqcup \cL_{-1,n}$,
the existence of $H: [0,1] \times K \to \cL_{+1} \sqcup \cL_{-1,n}$,
$H(0,\cdot) = f$, $H(1,\cdot) = \nu_2 \ast f$
or $H(1,\cdot) = F_{2n} \circ f$ is not clear at this point.
Indeed, the existence (or not) of such a homotopy
is the crucial point in this paper.
The following proposition helps clarify the situation.

\begin{prop}[\cite{S1}]
\label{prop:S1}
Let $K$ be a compact set and 
let $f: K \to \cL_I \subset \cI_I$ a continuous function.
\begin{enumerate}[(a)]
\item{For sufficiently large $n$,
the functions $\nu_2 \ast f$ and $F_{2n} \circ f$ are homotopic in $\cL_I$.}
\item{If $f = \nu_1 \ast \tilde f$ (for some $\tilde f$)
then $f$ is homotopic to $\nu_2 \ast f$.}
\item{The function $f$ is homotopic to a constant in $\cI_I$
if and only if $\nu_2 \ast f$ is homotopic to a constant in $\cL_I$.}
\item{There exists a map $\bfg_2: \Ss^2 \to \cL_{+1}$
such that $\nu_2 \ast \bfg_2$ is homotopic to a constant in $\cL_{+1}$
but $\bfg_2$ is not.}
\end{enumerate}
\end{prop}

\section{Bruhat cells and the set $\cW$}

Given a locally convex curve $\gamma: [0,1] \to \Ss^2$
and a $3\times 3$ matrix $A$ with positive determinant,
the curve $\gamma^A: [0,1] \to \Ss^2$, 
\[ \gamma^A(t) = \frac{A \gamma(t)}{|A \gamma(t)|} \]
is also locally convex. 
Furthermore, $\fF_{\gamma^A}(t) = A \fF_\gamma(t) U$
for $U \in \cUp$, where $\cUp$ is the group of upper triangular
$3\times 3$ matrices with positive off-diagonal entries.

Let $\cUpu \subset \cUp$ be the group of upper triangular matrices
with unit diagonal.
Recall that $SO(3)$ is divided in {\it Bruhat cells} by the following
equivalence relation: $Q_1$ and $Q_2$ are equivalent
if and only if there exist $U_1 \in \cUpu$ and $U_2 \in \cUp$
with $Q_1 = U_1 Q_2 U_2$.
The group Weyl group $D_3 \subset SO(3)$ of signed permutation matrices
with positive determinant has one element per cell.
The four open cells $\cJ_1, \cJ_2, \cJ_3, \cJ_4$ have respective
representatives
\begin{gather*}
J_1 = \begin{pmatrix} 0 & 0 & -1 \\ 0 & 1 & 0 \\ 1 & 0 & 0 \end{pmatrix}, \quad
J_2 = \begin{pmatrix} 0 & 0 & 1 \\ 0 & -1 & 0 \\ 1 & 0 & 0 \end{pmatrix}, \\
J_3 = \begin{pmatrix} 0 & 0 & 1 \\ 0 & 1 & 0 \\ -1 & 0 & 0 \end{pmatrix}, \quad
J_4 = \begin{pmatrix} 0 & 0 & -1 \\ 0 & -1 & 0 \\ -1 & 0 & 0 \end{pmatrix}.
\end{gather*}
Further recall (\cite{SK}, \cite{ShapiroM}, \cite{SaSha})
that given $Q \in SO(3)$ there exists a convex curve $\gamma: [0, 1]$
with $\fF_\gamma(0) = I$ and $\fF_\gamma(1) = Q$
if and only if $Q$ belongs to $\cJ_2$
or to one of the $5$ lower dimensional cells
in its boundary corresponding to the following matrices:
\[
\begin{pmatrix} 1 & 0 & 0 \\ 0 & 1 & 0 \\ 0 & 0 & 1 \end{pmatrix}, \quad
\begin{pmatrix} 0 & 1 & 0 \\ -1 & 0 & 0 \\ 0 & 0 & 1 \end{pmatrix}, \quad
\begin{pmatrix} 1 & 0 & 0 \\ 0 & 0 & 1 \\ 0 & -1 & 0 \end{pmatrix}, \quad
\begin{pmatrix} 0 & 0 & 1 \\ -1 & 0 & 0 \\ 0 & -1 & 0 \end{pmatrix}, \quad
\begin{pmatrix} 0 & 1 & 0 \\ 0 & 0 & 1 \\ 1 & 0 & 0 \end{pmatrix}. 
\]
As a consequence, given $Q \in SO(3)$, there exists $\gamma \in \cL_{-1,c}$
with $\fF_\gamma(1/2) = Q$ if and only if $Q \in \cJ_2$.
Similarly, $Q_0^{-1} Q_1 \in \cJ_4$
if and only if there exists a convex curve $\gamma: [0,1] \to \Ss^2$
with $\gamma(0) = \gamma(1) = Q_0e_1$, $\gamma'(0) = \gamma'(1) = -Q_0e_2$,
$\gamma(1/2) = Q_1e_1$, $\gamma'(1/2) = -Q_1e_2$.

Let $\cW \subset \cL_I \times \Ss^1 \times \Ss^1$
be the set of triples $(\gamma,t_0,t_1)$
such that $(\fF_\gamma(t_0))^{-1} \fF_\gamma(t_1) \in \cJ_4$.
Define $U_1: \cW \to \cUpu$ and $U_2: \cW \to \cUp$
so that $(\fF_\gamma(t_0))^{-1} \fF_\gamma(t_1) = U_1 J_4 U_2$.
Alternatively, $(\gamma, t_0, t_1) \in \cW$ if and only if there exists
a convex curve $\alpha: \Ss^1 \to \Ss^2$,
$\alpha(0) = \gamma(t_0)$, $\alpha'(0) = -\gamma'(t_0)$,
$\alpha(1) = \gamma(t_1)$, $\alpha'(1) = -\gamma'(t_1)$.
Notice that $\cW$ is an open subset of $\cL_I \times \Ss^1 \times \Ss^1$;
$(\gamma,t_0,t_1) \in \cW$ implies $t_0 \ne t_1$ and
$(\gamma,t_1,t_0) \in \cW$.

We now define the function 
$\Delta^\sharp: [0,+\infty) \times \cW \to \cL_I$,
one of our main technical tools throughout the paper.

\begin{defin}
\label{defin:Deltasharp}
Let $(\gamma, t_0, t_1) \in \cW$ and $t_\base \in \Ss^1$
with $t_\base < t_0 < t_1 < t_\base + 1$.
Take $\epsilon^\sharp > 0$,
$\epsilon^\sharp < (1/20) \min(t_0-t_\base,t_1-t_0,t_\base+1-t_1)$.
Let $Q_0 = \fF_\gamma(t_0)$, $Q_1 = \fF_\gamma(t_1)$,
$U_1 = U_1(\gamma,t_0,t_1)$, $U_2 = U_2(\gamma,t_0,t_1)$
so that $Q_0^{-1} Q_1 = U_1 J_4 U_2$.
Let $\alpha = \gamma^{U_1^{-1} Q_0^{-1}}$, i.e.,
\[ \alpha(t) =
\frac{U_1^{-1} Q_0^{-1} \gamma(t)}{|U_1^{-1} Q_0^{-1} \gamma(t)|}; \]
$\alpha$ is locally convex with
$\fF_\alpha(t_0) = I$, $\fF_\alpha(t_1) = J_4$.
For $s \in [0,+\infty)$, let $T_s = \min(s,2)$;
define an increasing piecewise linear homeomorphisms
$h_{a,s}: [0, T_s \epsilon^\sharp] \to [0,s]$
whose graph is a polygonal line
with vertices $(0,0)$, $(\epsilon^\sharp,1)$ (if $s \ge 1$) and
$(T_s \epsilon^\sharp,s)$.
Similarly, let the graph of
$h_{b,s}: [t_0 + \epsilon^\sharp T_s, t_1 - \epsilon^\sharp T_s]
\to [t_0,t_1]$
have vertices
$(t_0 + T_s \epsilon^\sharp, t_0)$,
$(t_0 + 4\epsilon^\sharp, t_0 + 4\epsilon^\sharp)$
$(t_1 - 4\epsilon^\sharp, t_1 - 4\epsilon^\sharp)$
and $(t_1 - T_s \epsilon^\sharp, t_1)$.
Define $\alpha_s: [t_0,t_1] \to \Ss^2$ by
\[ \alpha_s(t) = \begin{cases}
\nu_1\left(h_{a,s}(t-t_0)\right),&
t_0 \le t \le t_0 + T_s \epsilon^\sharp, \\
\fF_{\nu_1}(s) \alpha\left(h_{b,s}(t)\right),&
t_0 + T_s \epsilon^\sharp \le t \le t_1 - T_s \epsilon^\sharp, \\
J_4 \nu_1\left(- h_{a,s}(t_1-t)\right),&
t_1 - T_s \epsilon^\sharp \le t \le t_1. \end{cases} \]
Finally, 
\[ Q \Delta^\sharp(s,\gamma,t_0,t_1) = \begin{cases}
\gamma(t),& t_\base \le t \le t_0 \textrm{ or } t_1 \le t \le t_\base,\\
(\alpha_s)^{Q_0 U_1}(t),& t_0 \le t \le t_1, \end{cases} \]
where $Q \in SO(3)$ is uniquely chosen so that
$\fF_{\Delta^\sharp(s,\gamma,t_0,t_1)}(0) = I$.
\end{defin}

A few remarks are in order.
The curve $\alpha_s$ is obtained from $\alpha$ by attaching 
an arc of circle of angle $2\pi s$ to either end,
rotating the curve to keep the same endpoints and reparametrizing.
Similarly, $\Delta^\sharp(s,\gamma,t_0,t_1)$
is obtained from $\gamma$ by inserting
an arc of $2\pi s$ at positions $t_0$ and $t_1$;
if $s$ is not an integer, the portion of $\gamma$
between $t_0$ and $t_1$ will be ``rotated''.
Up to minor deformations, the path $\Delta^\sharp(s,\gamma,t_0,t_1)$,
$s \in [0,1]$, from $\gamma$ to $\gamma^\sharp_{t_0,t_1} =
\Delta^\sharp(1,\gamma,t_0,t_1)$ is exemplified
in Figure \ref{fig:little} (in a situation where $0 < t_0 < t_1 < 1$).
The loop between $t_0$ and $t_1$
is pushed along a geodesic all the way, until it comes back (a).
The two long chunks of curve (in the figure, very nearly geodesics)
are then shrunk (b) and rounded (c) so that we obtain $\gamma^\sharp_{t_0,t_1}$.
Notice that the portion of $\gamma$ outside the interval $[t_0,t_1]$
is unaltered throughout the process.

Strictly speaking, the definition of $\Delta^\sharp$
depends on the choice of $\epsilon^\sharp$:
the only difference, however, when you change $\epsilon^\sharp$
is that functions get reparametrized.
Similarly, we ask that $U_1 \in \cUpu$
so that $U_1$ and $U_2$ become uniquely determined.
The function $h_s$ in Definition \ref{defin:Deltasharp}
is chosen so that the following technical result holds.

\begin{lemma}
\label{lemma:compactcW}
Let $K$ be a compact set and $(f,t_0,t_1): K \to \cW$ a continuous map;
there exists $\epsilon^\sharp > 0$ which suits
the definition of $\Delta^\sharp(s,f(p),t_0(p),t_1(p))$ for all $s$.

Let $\cK \subset \cL_I$ be the compact set of all curves
of the form $\Delta^\sharp(s,f(p),t_0(p),t_1(p))$, $p \in K$, $s \in [0,3]$.
For $s > 0$, let $n = \lfloor s-3 \rfloor$, $\tilde s = s - n$.
For $p \in K$, let $\gamma = \Delta^\sharp(s,f(p),t_0(p),t_1(p))$,
$\tilde\gamma = \Delta^\sharp(\tilde s,f(p),t_0(p),t_1(p))$.
Let
\[ h(t) = \begin{cases} 0,& t_\base \le t \le t_0,\\
n(t-t_0)/{\epsilon^\sharp},& t_0 \le t \le t_0 + \epsilon^\sharp,\\
n,& t_0 + \epsilon^\sharp \le t \le t_1 - \epsilon^\sharp, \\
2n - n(t_1-t)/{\epsilon^\sharp},& t_1 - \epsilon^\sharp \le t \le t_1,\\
2n,& t_1 \le t \le t_\base. \end{cases} \]
Then $\tilde\gamma \in \cK$ and
\[ \gamma(t) =
Q \fF_{\tilde\gamma}(t) \nu^{U_1(p)}_1(h(t)), \quad Q \in SO(3), \]
where $U_1(p) = U_1(f(p),t_0(p),t_1(p))$.
\end{lemma}

{\nobf Proof: } This is a straightforward computation. \qed

The functions
\( \Delta^\sharp(\tilde s,(\Delta^\sharp(s,\gamma,t_0,t_1),t_0,t_1)) \)
and \( \Delta^\sharp(s+\tilde s,\gamma,t_0,t_1) \)
differ by reparametrization only.
Under suitable hypothesis, a related identity holds 
for distinct points $(\gamma,t_0,t_1)$, $(\gamma,t_2,t_3)$.
Two points $(\gamma,t_0,t_1), (\gamma, t_2, t_3) \in \cW$ are
\textit{disjoint} if the intervals $[t_0,t_1], [t_2,t_3] \subset \Ss^1$
are disjoint, or, equivalently, if
$t_0 < t_1 < t_2 < t_3 < t_0+1$ or $t_2 < t_3 < t_0 < t_1 < t_2+1$.

\begin{lemma}
\label{lemma:disjDeltasharp}
If $(\gamma,t_0,t_1), (\gamma,t_2,t_3) \in \cW$ are disjoint then,
for any $s, \tilde s \in [0,1]$,
\[ (\Delta^\sharp(s,\gamma,t_0,t_1),t_2,t_3),
(\Delta^\sharp(\tilde s,\gamma,t_2,t_3),t_0,t_1) \in \cW. \]
Furthermore,
\[ \Delta^\sharp(\tilde s,(\Delta^\sharp(s,\gamma,t_0,t_1),t_2,t_3)) =
\Delta^\sharp(s, (\Delta^\sharp(\tilde s,\gamma,t_2,t_3),t_0,t_1)). \]
\end{lemma}

{\nobf Proof: }
This follows directly from the construction of $\Delta^\sharp$
and of the function $\gamma^\ast_s$ in the definition.
Indeed, $\gamma^\ast_s$ for $(t_0,t_1)$ coincides with $\gamma$
in an open interval containing $(t_2,t_3)$ (and vice versa).
\qed

\section{Disjoint covers}

Recall that one of our aims is to decide whether
a continuous map $f: K \to \cL_I$ is homotopic to $\nu_2 \ast f$.
In this section, we present several situations where this is the case
and one example where this is not the case.
We start with a simple example.

\begin{lemma}
\label{lemma:cW0}
Let $K$ be a compact set and $f: K \to \cL_I$ a continuous map.
If there is a continuous function $t_1: K \to \Ss^1$ such that
for all $p \in K$ we have $(f(p),0,t_1(p)) \in \cW$
then $f$ is homotopic to $\nu_2 \ast f$.
\end{lemma}

{\nobf Proof: }
Let $H: [0,1] \times K \to \cL_I$ be defined by
$H(s,p) = \Delta^\sharp(2s,f(p),0,t_1(p))$.
Up to reparametrization, $H(1,p) = \nu_2 \ast \tilde f(p)$,
where $\tilde f(p)$ is obtained from $f(p)$ by inserting two turns at $t_1(p)$.
Since $\nu_2$ and $\nu_4$ are in the same connected component,
$\nu_2 \ast \tilde f$ is homotopic to $\nu_2 \ast (\nu_2 \ast \tilde f)$
and therefore to $\nu_2 \ast f$.
\qed

Given this result, some questions are natural:
\begin{itemize}
\item{If there exist continuous functions $t_0, t_1: K \to \Ss^1$
with $(f(p),t_0(p),t_1(p)) \in \cW$, does it follow that
$f$ is homotopic to $\nu_2 \ast f$?}
\item{If for every $p$ there exist $t_0, t_1 \in \Ss^1$
such that $(f(p),t_0,t_1) \in \cW$, does it follow that
$f$ is homotopic to $\nu_2 \ast f$?}
\end{itemize}
As we shall see, the answers are yes and no, respectively.
Before we attack these problems, however, we introduce a few concepts.

Define 
\[ \cO = \{ \gamma \in \cL_I \;|\; \forall t_0, t_1 \in \Ss^1,
(\gamma, t_0, t_1) \notin \cW \}: \]
the set $\cO$ is clearly a closed subset of $\cL_I$.
A {\it double point} of a curve $\gamma \in \cL_I$
is a pair $(t_0, t_1) \in (\Ss^1)^2$, $t_0 \ne t_1$,
with $\gamma(t_0) = \gamma(t_1)$.
Similarly, an {\it $n$-tuple point} is an $n$-tuple
$(t_0, t_1, \ldots, t_{n-1})$,
$0 \le t_0 < t_1 < \cdots < t_{n-1} < 1$, such that
$\gamma(t_0) = \gamma(t_1) = \cdots = \gamma(t_{n-1})$.
We identify the double point $(t_1, t_0)$ with $(t_0,t_1)$.
A double point $(t_0, t_1)$ is  a {\it self-tangency} if
$\gamma'(t_0)$ and $\gamma'(t_1)$ are parallel
and {\it transversal} otherwise.
A self-tangency $(t_0, t_1)$ is {\it positive} if $\gamma'(t_1)$
is a positive multiple of $\gamma'(t_0)$ (and {\it negative} otherwise).

If $(t_0,t_1)$ is a positive self-tangency, the normal vector
$\bfn_\gamma(t) = \fF_\gamma(t) e_3$
satisfies $\bfn_\gamma(t_0) = \bfn_\gamma(t_1)$.
Define $h_1(t) = \langle \gamma(t), \fF_\gamma(t) e_2 \rangle$,
$h_2(t) = \langle \gamma(t), \bfn_\gamma(t_0) \rangle$:
notice that $h_1(t_0) = h_1(t_1) = h_2(t_0) = h_2(t_1) = 0$,
$h_2'(t_0) = h_2'(t_1) = 0$, $h_1'(t_0) > 0$ and $h_1'(t_1) > 0$.
There exists $\epsilon > 0$ such that
$h_1|_{(t_0-\epsilon, t_0+\epsilon)}$ and 
$h_1|_{(t_1-\epsilon, t_1+\epsilon)}$ are invertible.
Let $g_0 = h_2 \circ (h_1|_{(t_0-\epsilon, t_0+\epsilon)})^{-1}$ and 
$g_1 = h_2 \circ (h_1|_{(t_1-\epsilon, t_1+\epsilon)})^{-1}$:
near the origin, the graph of $g_i$ is the orthogonal projection of
the image of the curve $\gamma$ near $\gamma(t_i)$ to the tangent plane,
using $\gamma'(t_i)/|\gamma'(t_i)|$ and $\bfn_\gamma(t_i)$ as basis.
By construction, $g_0(0) = g_1(0) = 0$ and $g_0'(0) = g_1'(0) = 0$.
By convexity, $g_0''(0) > 0$, $g_1''(0) > 0$.
The self-tangency $(t_0,t_1)$ has {\it order $n$}
if $g_0^{(n)}(0) \ne g_1^{(n)}(0)$
but $g_0^{(j)}(0) = g_1^{(j)}(0)$ for any $j < n$
and order $+\infty$ if $g_0^{(j)}(0) = g_1^{(j)}(0)$ for any $j$.
A {\it self-osculating} point is a positive self-tangency of order $3$ or more.

\begin{lemma}
\label{lemma:arctoloop}
Let $(t_0, t_1)$ be a double point of $\gamma \in \cL_I$.
If $(t_0,t_1)$ is either transversal or a negative self-tangency
then for any $\epsilon > 0$ there exist $\tilde t_0, \tilde t_1$,
$|\tilde t_0 - t_0| + |\tilde t_1 - t_1| < \epsilon$,
$(\gamma,t_0,t_1) \in \cW$.
\end{lemma}

\begin{figure}[ht]
\begin{center}
\psfrag{t0-}{$t_0 - \epsilon/4$}
\psfrag{t0+}{$t_0 + \epsilon/4$}
\psfrag{t1-}{$t_1 - \epsilon/4$}
\psfrag{t1+}{$t_1 + \epsilon/4$}
\epsfig{height=30mm,file=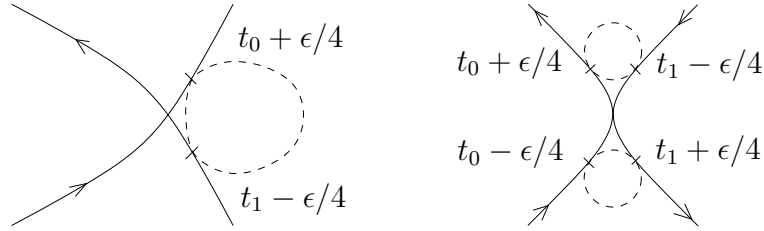}
\end{center}
\caption{Obtaining $(\tilde t_0, \tilde t_1)$ such that
$(\gamma, \tilde t_0, \tilde t_1) \in \cW$.}
\label{fig:arctoloop}
\end{figure}

{\nobf Proof: }
First consider transversal double points.
Assume without loss of generality that
$\det(\gamma(t_0), \gamma'(t_0), \gamma'(t_1)) > 0$.
For small $\epsilon$ we may take
$\tilde t_0 = t_0 + \epsilon/4$, $\tilde t_1 = t_1 - \epsilon/4$:
this is a straightforward computation but is probably
best verified geometrically in Figure \ref{fig:arctoloop}:
the dashed convex curves tangent to $\gamma$ validate
the geometric characterization of $\cW$.

For a negative self-tangency and small $\epsilon$ we may take either
$\tilde t_0 = t_0 + \epsilon/4$, $\tilde t_1 = t_1 - \epsilon/4$ or
$\tilde t_0 = t_0 - \epsilon/4$, $\tilde t_1 = t_1 + \epsilon/4$:
this is again verified in Figure \ref{fig:arctoloop}.
\qed

It follows directly from this result that if $\gamma \in \cO$
then all double points of $\gamma$ are positive self-tangencies.

\begin{figure}[ht]
\begin{center}
\psfrag{t0}{$t_0$}
\epsfig{height=25mm,file=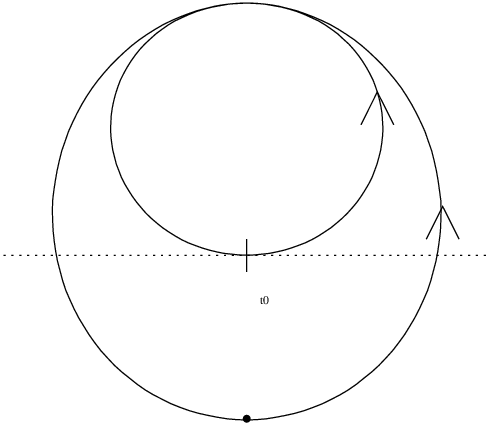}
\end{center}
\caption{A curve in $\cO_2$.}
\label{fig:oscobruhat}
\end{figure}

The following proposition settles the second question raised
at the beginning of this section.
Recall (\cite{S1}) that $\bfg_{+,2}: \Ss^2 \to \cL_{+1}$
is an explicit function such that $\bfg_{+,2}$ and $\nu_2 \ast \bfg_{+,2}$
are not homotopic in $\cL_{+1}$
(even though they are homotopic in $\cI_{+1}$).

\begin{prop}
\label{prop:Oisnotenough}
There is a map $f: \Ss^2 \to \cL_{+1} \smallsetminus \cO$
which is homotopic (in $\cL_{+1}$) to $\bfg_2$.
In particular, $f$ and $\nu_2 \ast f$ are not homotopic.
\end{prop}

\begin{figure}[ht]
\begin{center}
\epsfig{height=120mm,file=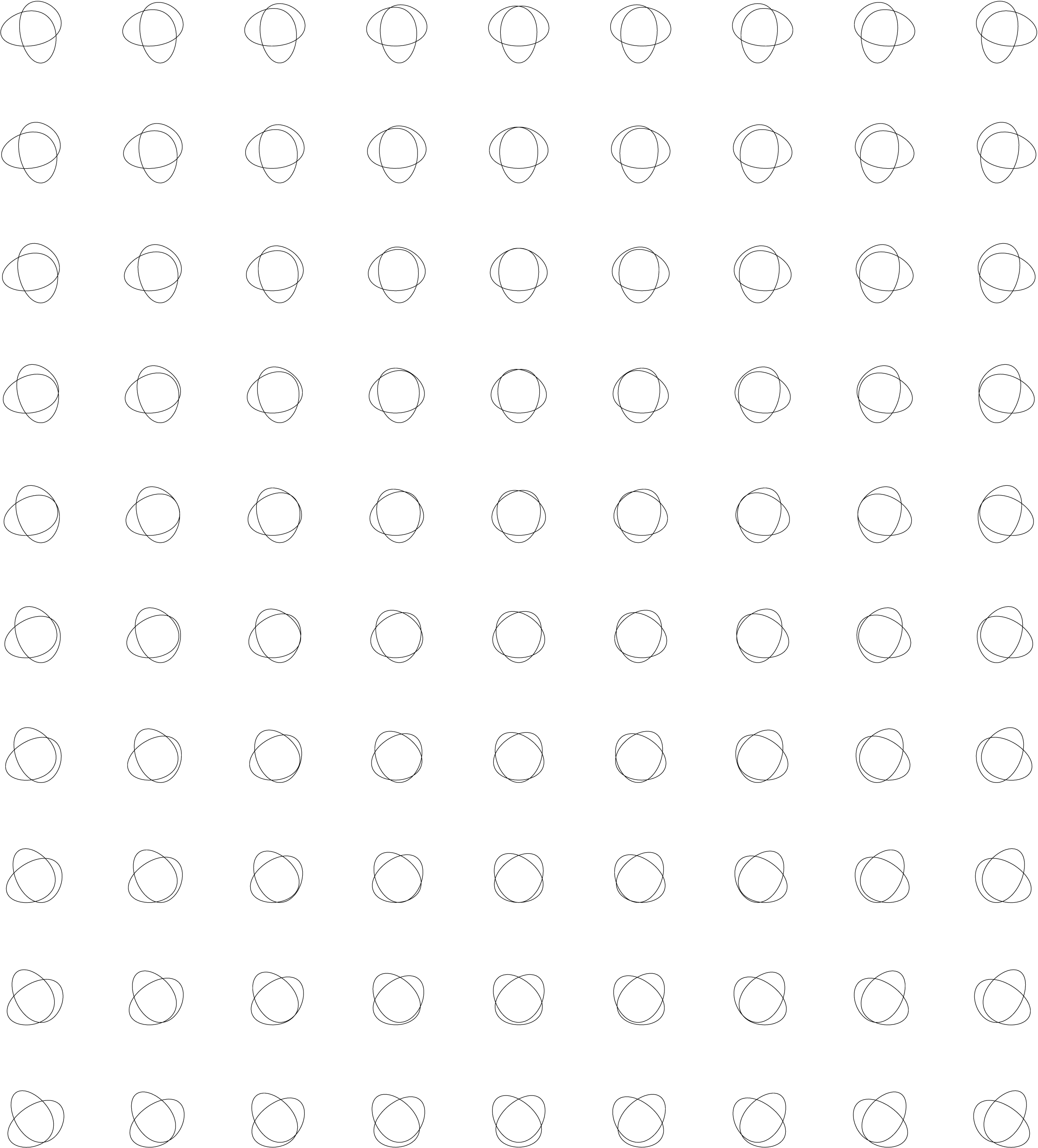}
\end{center}
\caption{A family of convex curves.}
\label{fig:pentaB}
\end{figure}

{\nobf Proof: }
Let $\Ss^1$ be the unit circle in the complex plane
and let $D \subset \CC$ be the closed disk of radius $1/4$.
For $s \in D$, let $g_s: \Ss^1 \to \RR^2$ be defined by
\[ g_s(z) = -i z^2 - \frac{i}{10} z^{-3} + s z^{-1}; \]
some such curves are drawn in Figure \ref{fig:pentaB}
(the curve corresponding to $s = 0$ is in position $(6,5)$,
i.e., sixth row, fifth column).

A straightforward computation verifies that these curves are locally convex
in the plane: central projection obtains a similar family
of locally convex curves in the sphere.
Notice that there are $5$ curves with self-osculating points
(approximately in positions $(3,5)$, $(5,2)$, $(5,8)$, $(9,3)$ and $(9,7)$)
but none of them lie in $\cO$ since they all have transversal double points.

Remove a small disk near $\nu_2$ in the function $\bfg_{+,2}$
corresponding to the region below the bottom row in Figure 9 of \cite{S1};
the above family can be used to plug the hole.
The resulting function if homotopic to $\bfg_{+,2}$
and its image is contained in $\cL_{+1} \smallsetminus \cO$.
\qed

Before we answer the first question,
we present another situation where $f$ is guaranteed
to be homotopic to $\nu_2 \ast f$.

\begin{lemma}
\label{lemma:almostnu1}
Let $K$ be a compact manifold and $f: K \to \cL_I$ be a continuous map.
Assume there exist functions $t_0, t_1: K \to (0,1)$ such that,
for all $p \in K$:
\begin{enumerate}[(a)]
\item{$0 < t_0(p) < t_1(p) < 1$;}
\item{$f(p)|_{[0,t_0(p)]}$ and $f(p)|_{[t_0(p),t_1(p)]}$ are convex;}
\item{there exists a convex curve $\alpha(p): [0,1] \to \Ss^2$
with $\fF_{\alpha(p)}(0) = \fF_{f(p)}(t_0(p))$,
$\fF_{\alpha(p)}(1/2) = \fF_{f(p)}(0)$,
$\fF_{\alpha(p)}(1) = \fF_{f(p)}(t_1(p))$.}
\end{enumerate}
Then $f$ is homotopic to $\nu_2 \ast f$.
\end{lemma}

{\nobf Proof: }
Define $H: [0,1] \to \cL_I$ for $s \in [0,1/2]$ with $H(0,p) = f(p)$ and
\[ H(1/2,p)(t) = \begin{cases} f(p)(t), & t \in [0,t_0(p)] \cup [t_1(p),1],\\
\alpha(p)\left( \frac{t-t_0(p)}{t_1(p)-t_0(p)} \right),& t \in [t_0(p),t_1(p)].
\end{cases} \]
Contractibility of the space of convex curves with prescribed
initial and final value and direction guarantees that this can be done.
Let $t_{1/2}(0) = (t_0(p) + t_1(p))/2$:
notice that $H(1/2,p)|_{[0,t_{1/2}(p)]}$ is a closed convex curve.
Set
\[ H(1,p)(t) = \begin{cases} \nu_1(2t), & t \in [0,1/2], \\
H(1/2,p)\left(\frac{t - t_{1/2}(p)}{1 - t_{1/2}(p)}\right),& t\in[1/2,1];
\end{cases} \]
contractibility of $\cL_{-1,c}$ guarantees that this can be done.
Now $\tilde f: K \to \cL_I$, $\tilde f(p) = H(1,p)$,
is of the form $\tilde f = \nu_1 \ast \textrm{(something)}$.
The result now follows from Proposition \ref{prop:S1}, item (b).
\qed

For a function $f: K \to \cL_I$,
a finite open cover $K = \bigcup_{i=1, \ldots, N} V_i$
together with functions $t_{0,i}, t_{1,i}: V_i \to \Ss^1$
is a \textit{disjoint cover of $f$} if:
\begin{enumerate}[(a)]
\item{if $p \in V_i$ then $(f(p),t_{0,i}(p), t_{1,i}(p)) \in \cW$;}
\item{if $p \in V_i \cap V_j$ then the points
$(f(p),t_{0,i}(p),t_{1,i}(p)), (f(p),t_{0,j}(p),t_{1,j}(p)) \in \cW$
are either equal or disjoint.}
\end{enumerate}

The following lemma settles the fist question;
it is actually much stronger.

\begin{lemma}
\label{lemma:disjcover}
Let $K$ be a compact manifold and $f: K \to \cL_I$ be a continuous map.
If the function $f$ admits a disjoint cover
then $f$ is homotopic to $\nu_2 \ast f$.
\end{lemma}

It is instructive to verify directly that the function $f$ constructed
in Proposition \ref{prop:Oisnotenough} does not admit a disjoint cover.

{\nobf Proof: }
Consider a function $f$ and a cover by disjoint loops with the notation above.
Recall that the support $\supp(\phi)$ of $\phi: K \to [0,1]$
is the closure of $\phi^{-1}((0,1])$.
Let $\phi_i: K \to [0,1]$ be a partition of unit:
$\supp(\phi_i) \subset V_i$ and $\sum_i \phi_i(p) = 1$.

Intuitively, our first step is to apply $\Delta^\sharp(\sigma,\cdots)$
to the points given by the cover;
$\sigma$ is a function which goes to zero together with $\phi_i$.
More precisely,
select $\epsilon^\sharp > 0$ such that if $\phi_{i_1}(p), \phi_{i_2}(p) > 0$
and $t_{0,i_1}(p) < t_{1,i_1}(p) <
t_{0,i_2}(p) < t_{1,i_2}(p) < 1 + t_{0,i_1}(p)$ then
\[ 20 \epsilon^\sharp < \min(t_{1,i_1}(p) - t_{0,i_1}(p),
t_{0,i_2}(p) - t_{1,i_1}(p), t_{1,i_2}(p) - t_{0,i_2}(p),
1 + t_{0,i_1}(p) - t_{1,i_2}(p)). \]
For $p \in K$, let $(t_{0,j,p},t_{1,j,p})$,
$j = 1, \ldots, N_p$, $1 \le N_p \le N$,
be the distinct pairs for which there exists $i \in \{1, \ldots, N\}$,
$p \in V_i$, $t_{0,j,p} = \tilde t_{0,i}(p)$, $t_{1,j,p} = \tilde t_{1,i}(p)$. 
For each $j$, let $I_j \subset \{1, \ldots, N\}$ be the set of indices $i$
for which the above conditions hold.
Define \[ \psi_j(p) = \sum_{i \in I_j} \phi_i(p). \]
For $s \in [0,1/2]$, let $\sigma_j(s,p) = \min(2M,12MNs \psi_j(p))$
where $M > 0$ is a large integer to be specified later.
Let $H_0(s,p) = f(p)$ and define recursively
\[H_j(s,p) = \Delta^\sharp(\sigma_j(s,p),H_{j-1}(s,p),t_{0,j,p},t_{1,j,p})\]
and $H(s,p) = H_{N_p}(s,p)$.
By Lemma \ref{lemma:disjDeltasharp} the order of the indices does not matter.

Intuitively, we added many turns to each curve and must now spread them.
Define $U_1: K \times \Ss^1 \to \cUpu$ so that if $p \in V_i$,
$t \in [t_{0,i}(p), t_{0,i}(p) + 4 \epsilon^\sharp] \cup
[t_{1,i}(p) - 4 \epsilon^\sharp, t_{1,i}(p)]$ then
$U_1(p,t) = U_1(f(p),t_{0,i}(p),t_{1,i}(p))$.
Write
\[ R(p,t,x) = U_1(p,t) \fF_{\nu_1}(x) (U_1(p,t))^{-1}; \]
notice that $R(p,t,x_1) R(p,t,x_2) = R(p,t,x_1+x_2)$.
For any $p \in K$ there is at least one index $j$ such that
$\psi_j(p) \ge 1/(3N)$ and therefore $\sigma_i(1/2,p) = 2M$:
assume without loss of generality that these indices
are $j = 1, \ldots, \tilde N_p$, $1 \le \tilde N_p \le N_p$.
We define an auxiliary curve $\eta(p) \in \cL_I$ by
$\eta_0 = f(p)$,
$\eta_j = \Delta^\sharp(M,\eta_{j-1},t_{0,j,p},t_{1,j,p})$ for
$j = 1, \ldots, \tilde N_p$ and
$\eta_j = \Delta^\sharp(\sigma_j(1/2,p),\eta_{j-1},t_{0,j,p},t_{1,j,p})$ for
$j = \tilde N_p + 1, \ldots, N_p$ and $\eta(p) = \eta_{N_p}$ so that,
by Lemma \ref{lemma:compactcW},
\[ H(1/2,p)(t) = r(p,t) \fF_{\eta(p)}(t)
R\left(p,t,
M \sum_{j=1,\ldots,\tilde N_p} ( \beta_{0,j}(t) + \beta_{1,j}(t) ) \right)
e_1,
\]
where $r(p,t) \in (0,+\infty)$ and
\[
\beta_{0,j} = \beta_{[t_{0,j,p},t_{0,j,p}+\epsilon^\sharp]}, \quad
\beta_{1,j} = \beta_{[t_{1,j,p}-\epsilon^\sharp,t_{1,j,p}]}, \quad
\beta_{[t_-,t_+]}(t) = \begin{cases} 0, & t \le t_-, \\
\frac{t-t_-}{t_+-t_-}, & t_- \le t \le t_+, \\
1, & t \ge t_+. \end{cases} \]
Given an interval $[t_- , t_+] \subset \Ss^1$
define $t_\base = (t_+ + t_- - 1)/2$
so that $t_\base < t_- < t_+ < t_\base + 1$.
Given $\theta \in [0,1]$, let
\[ [t_-, t_+]^\theta
= [ \theta t_\base + (1-\theta) t_-,\theta (t_\base+1) + (1-\theta) t_+] \]
so that if $\theta = 1$ the interval degenerates to the whole circle.
Let $\theta_j: [1/2,1] \times K \to [0,1]$ be
\[ \theta_j(s,p) = \begin{cases}
0,& \psi_j(p) \le \frac{1}{3N},\\
(2s-1)(3N\psi_j(p) - 1),& \frac{1}{3N} \le \psi_j(p) \le \frac{2}{3N},\\
2s - 1,& \psi_j(p) \ge \frac{2}{3N}. \end{cases} \]
Define
\[ \beta_{0,j,s} =
\beta_{[t_{0,j,p},t_{0,j,p}+\epsilon^\sharp]^{\theta_j(s,p)}},
\quad
\beta_{1,j,s} =
\beta_{[t_{1,j,p}-\epsilon^\sharp,t_{1,j,p}]^{\theta_j(s,p)}},\]
\begin{align*}
H(s,p)(t) &= r(p,s,t) \fF_{\eta(p)}(t)
R\left(p,t,
M \sum_{j=1,\ldots,\tilde N_p} ( \beta_{0,j,s}(t) + \beta_{1,j,s}(t) ) \right)
e_1 \\
&= r(p,s,t) \fF_{\eta(p)}(t) U_1(p,t) \nu_1\left(
M \sum_{j=1,\ldots,\tilde N_p} ( \beta_{0,j,s}(t) + \beta_{1,j,s}(t) ) \right)
\end{align*}
where $r(p,s,t)$ is a positive number chosen so that
the expression has absolute value $1$.
For sufficiently large $M$, all the functions constructed above
will belong to $\cL_I$, as required.
Indeed, $\gamma = H(1/2,p)$ is of the form
$\gamma(t) = \gamma_1(t)/|\gamma_1(t)|$,
$\gamma_1(t) = \fF_\eta(t) U_1(p,t) \nu_1(h(t))$
where $\eta$ belongs to a compact set $\cK \subset \cL_I$
independent of the choice of $M$.
For $s \in [1/2,1]$, $\gamma = H(s,p)$ still has the same form;
for any given $t$, $s$ and $p$,
either $h(t) \in \ZZ$, $h'(t) = 0$ or
$h'(t) > M$ (up to a few transition points which need not concern us).
In the first case, local convexity of $\gamma$ follows from
local convexity of $H(1/2,p)$.
In the second case, expanding $\det(\gamma_1(t),\gamma_1'(t),\gamma_1''(t))$
shows that this expression is positive provided $M$ is large enough.

Let $\tilde f: K \to \cL_I$, $\tilde f(p) = H(1,p)$.
We claim that if $M$ is large enough then $\tilde f$
satisfies the hypothesis of Lemma \ref{lemma:almostnu1}.
Indeed, for $\gamma$ in the image of $\tilde f$, let
\[ \begin{pmatrix} g_1(t) \\ g_2(t) \\ g_3(t) \end{pmatrix}
= \begin{pmatrix} \sqrt{2}/2 & 0 & \sqrt{2}/2 \\
0 & 1 & 0 \\ -\sqrt{2}/2 & 0 & \sqrt{2}/2 \end{pmatrix}
(U_1(p,0))^{-1} (\fF_\gamma(0))^{-1} \gamma(t). \]
From the above form for $\gamma$ we have, for small $t$,
\[ \begin{pmatrix} g_1(t) \\ g_2(t) \\ g_3(t) \end{pmatrix} \approx
\frac{\sqrt{2}}{2} r(p,s,t)
\begin{pmatrix} \cos(2\pi h(t)) \\ \sin(2\pi h(t)) \\ 1 \end{pmatrix}. \]
Let $t_0(p)$ be the smallest $t>0$ for which $g_1(t) = 0$, $g_2(t) < 0$;
let $t_1(p)$ be the smallest $t > t_0(p)$
for which $g_1(t) = 0$, $g_2(t) > 0$.
For sufficiently large $M$ these are continuously defined
and the convexity hypothesis will hold.
This completes the proof of the claim and of the lemma.
\qed

\section{Loops and stars}

From now on our aim is to produce disjoint covers for functions
$f: K \to \cL_I$ or, at least, to prove that $f$ is homotopic to $\tilde f$
such that $\tilde f$ admits a disjoint cover.
We must therefore turn to the geometry of curves.

We now define nested dense open sets $\cL^{(k)}_{\pm 1} \subset \cL_{\pm 1}$
for $k \le 3$; the complement $\cL_{\pm 1} \smallsetminus \cL^{(k)}_{\pm 1}$
has codimension $k$.

\begin{figure}[ht]
\begin{center}
\epsfig{height=40mm,file=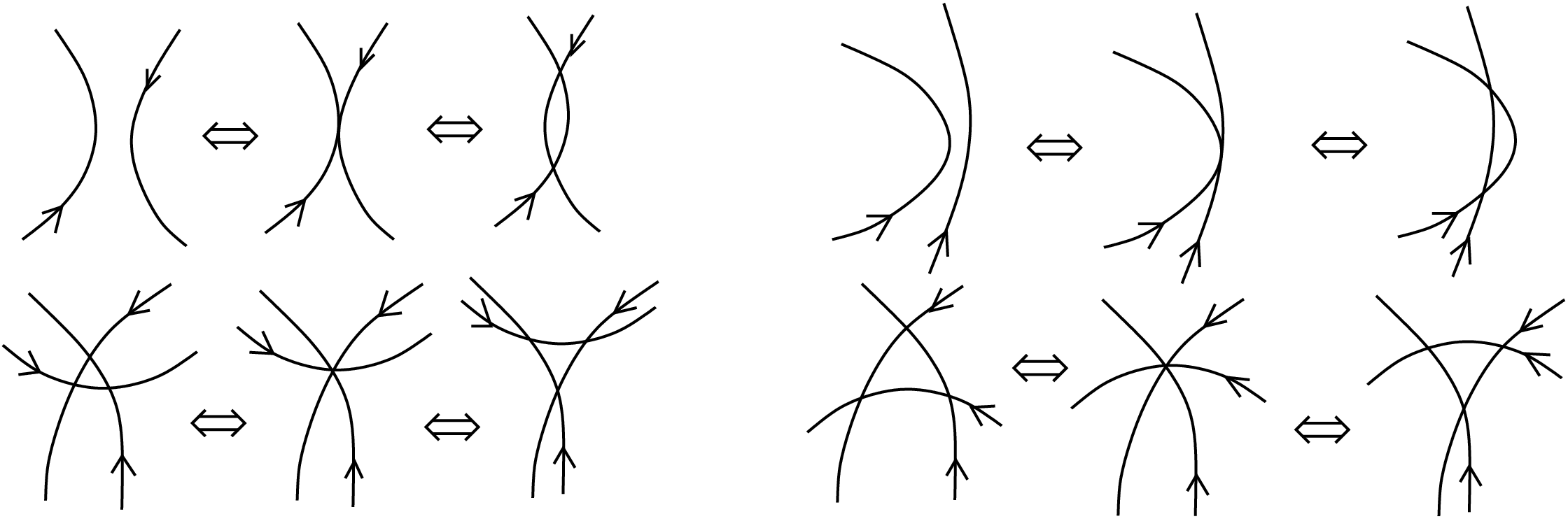}
\end{center}
\caption{Reidemeister moves: type II on first line, type III on second line.}
\label{fig:reidemeister}
\end{figure}

Let $\cL^{\langle 0 \rangle}_{\pm 1} = \cL^{(1)}_{\pm 1} \subset \cL_{\pm 1}$
be the set of curves with no triple points or self-tangencies.
Notice that the sets $\cL^{(1)}_{\pm 1}$ have infinitely 
many connected components
since the number of double points does not change in a connected
component of these sets.
Let $\cL^{\langle 1,a \rangle}_{\pm 1} \subset \cL_{\pm 1}$
be the set of curves with exactly one self-tangency of order $2$,
no triple points and
no self-osculating points (and an arbitrary number of double points).
Let $\cL^{\langle 1,b \rangle}_{\pm 1} \subset \cL_{\pm 1}$
be the set of curves with exactly one triple point and no self-tangencies.
The sets $\cL^{\langle 1, a \rangle}_{\pm 1},
\cL^{\langle 1, b \rangle}_{\pm 1} \subset \cL_{\pm 1}$
are disjoint submanifolds of codimension $1$.
Generically, the passage from one connected component
of $\cL^{(1)}_{\pm 1}$ to another
crosses $\cL^{\langle 1, \ast \rangle}_{\pm 1}$ transversally
and is a Reidemeister move
of type II (resp. III) if $\ast = a$ (resp. $\ast = b$; \cite{BZ});
Reidemeister moves of type I are not allowed in $\cL_I$.
Figure \ref{fig:reidemeister} shows the possible Reidemeister moves
in $\cL_I$.

\begin{figure}[ht]
\begin{center}
\psfrag{(b)}{$(b)$}
\psfrag{(c)}{$(c)$}
\psfrag{(d)}{$(d)$}
\psfrag{(e)}{$(e)$}
\epsfig{height=30mm,file=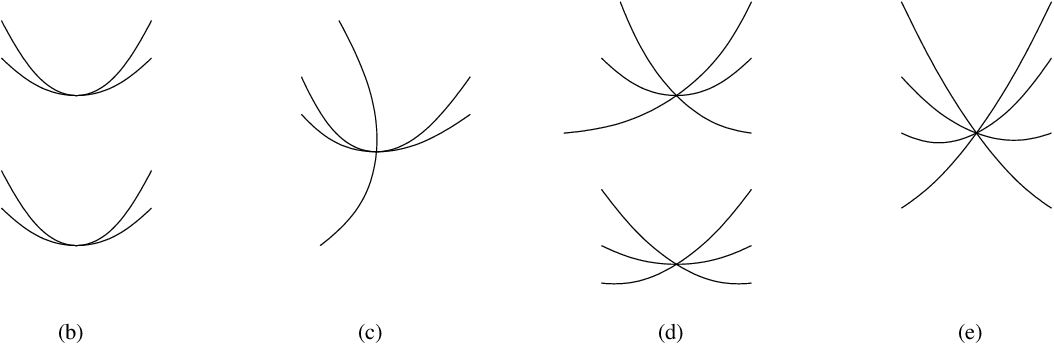}
\end{center}
\caption{Curves in $\cL^{\langle 2, \ast \rangle}_{\pm 1}$.}
\label{fig:cl2}
\end{figure}

Define $\cL^{(2)}_{\pm 1} = \cL^{(1)}_{\pm 1} \sqcup
\cL^{\langle 1 \rangle}_{\pm 1} \subset \cL_{\pm 1}$.
Let $\cL^{\langle 2, \ast \rangle}_{\pm 1}$ be the set of curves
having:
\begin{enumerate}[(a)]
\item{exactly one self-tangency, which is positive and of order $3$,
and no triple points;}
\item{exactly two self-tangencies, both of order $2$,
and no triple points;}
\item{exactly one triple point
where there is also a self-tangency of order $2$;}
\item{exactly two (unrelated) triple points;}
\item{exactly one quadruple point;}
\end{enumerate}
Figure \ref{fig:cl2} illustrates these situations.
These sets are submanifolds of codimension $2$.
Finally, define $\cL^{(3)}_{\pm 1} = \cL^{(2)}_{\pm 1} \sqcup
\cL^{\langle 2, a \rangle}_{\pm 1} \sqcup \cdots \sqcup
\cL^{\langle 2, e \rangle}_{\pm 1} \subset \cL_{\pm 1}$.

We already saw examples curves in $\cL^{\langle 2,a \rangle}_{+1}$
in Figures \ref{fig:oscobruhat} and \ref{fig:pentaB};
the second one shows a surface transversal to
$\cL^{\langle 2,a\rangle}_{+ 1}$.
Notice that the self-osculating point may be perturbed
to become one single transversal double point
or three transversal double points (observe the central column);
the reader should compare this with perturbations
of the real polynomial $P(x) = x^3$,
which may admit one or three real roots.

\begin{figure}[ht]
\begin{center}
\psfrag{(a)}{(a)}
\psfrag{(b)}{(b)}
\psfrag{t0}{$t_0$}
\psfrag{t1}{$t_1$}
\epsfig{height=25mm,file=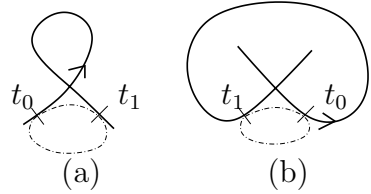}
\end{center}
\caption{A direct loop and a reverse loop.}
\label{fig:littlebig}
\end{figure}

A {\it loop} of a curve $\gamma \in \cL_I$
is a transversal double point $(t_0,t_1)$
such that the restriction $\gamma|_{[t_0,t_1)}$ is injective.
We sometimes think of the loop as the interval $[t_0,t_1]$,
the restriction of $\gamma$ to this interval
or even the image of this restriction.
A loop is {\it direct} (resp. {\it reverse}) if
$\det(\gamma(t_0), \gamma'(t_0), \gamma'(t_1))$ is negative (resp. positive).
Figure \ref{fig:littlebig} shows examples of loops.

\begin{figure}[ht]
\begin{center}
\epsfig{height=15mm,file=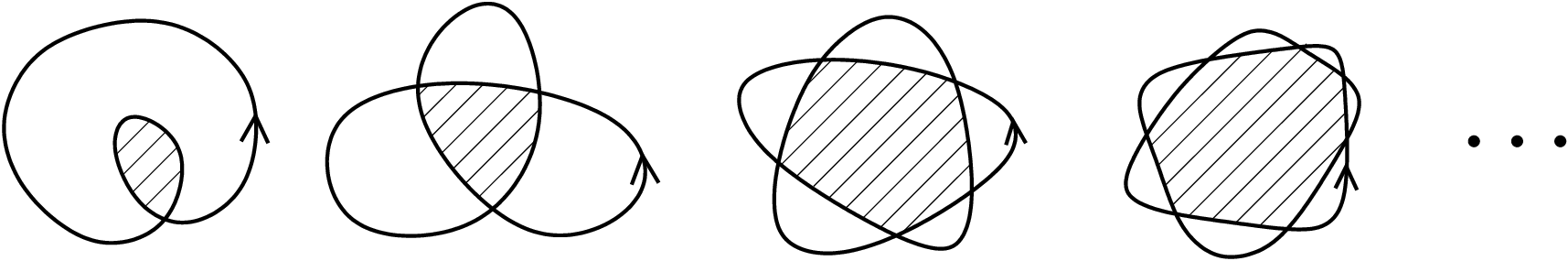}
\end{center}
\caption{Stars ($k = 0, 1, 2, 3, \ldots$).}
\label{fig:star}
\end{figure}

A {\sl star} is a curve $\gamma$ in the same connected component
of $\cL^{(1)}_{+1}$ as one of the infinite family of curves given in
Figure \ref{fig:star}.
More precisely, a star has $2k+1$ double points;
if $k > 0$, their images in the sphere are the vertices
of a convex polygon and, for any pair of adjacent vertices,
there are two arcs of $\gamma$ joining them.
Alternatively, a star is a curve in $\cL^{(1)}_{+1}$
which admits loops $(t_0,t_1)$ and $(t_1,t_0+1)$.

Let $\cT_0$ be the closure (in $\cL_{+1}$) of the set of stars
and let $\cT_1$ be its boundary.

A curve $\gamma \in \cL^{\langle 1,b \rangle}_{+1}$
with triple point $(t_0,t_1,t_2)$ is a {\sl trefoil} if
$(t_0,t_1)$, $(t_1,t_2)$ and $(t_2,t_0+1)$ are
direct loops.

\begin{lemma}
\label{lemma:trefoil}
The set $\cT_1$ is the set of trefoils and is a manifold of codimension 1.
The set $\cT_0$ is contractible and $\cT_1$ is homotopically equivalent
to $\Ss^1$.
\end{lemma}

\begin{figure}[ht]
\begin{center}
\epsfig{height=15mm,file=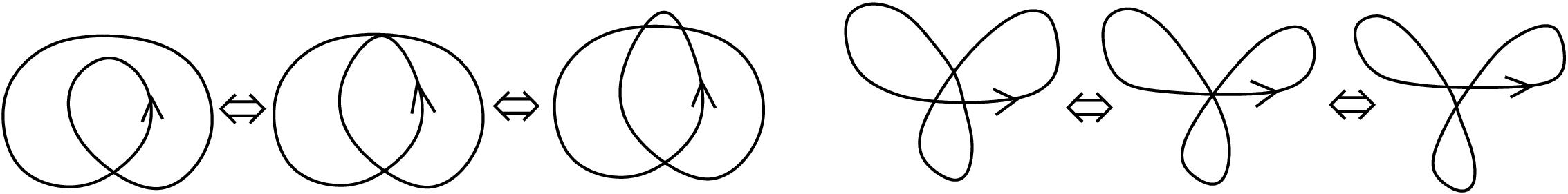}
\end{center}
\caption{Reidemeister moves starting at a star.}
\label{fig:reidstar}
\end{figure}

{\nobf Proof: }
We have to show that the only Reidemeister moves from a star
to a generic $\gamma$ which is not a star pass through a trefoil.
In order to do this, we classify all possible Reidemeister moves
starting at a star.
Figure \ref{fig:reidstar} shows how a Reidemeister move of type II
takes a star to another star (changing the value of $k$)
and how a Reidemeister move of type III takes a star ($k = 1$)
to a generic curve which is not a star passing through a trefoil.
We prove that these are the only possible moves.

\begin{figure}[ht]
\begin{center}
\epsfig{height=50mm,file=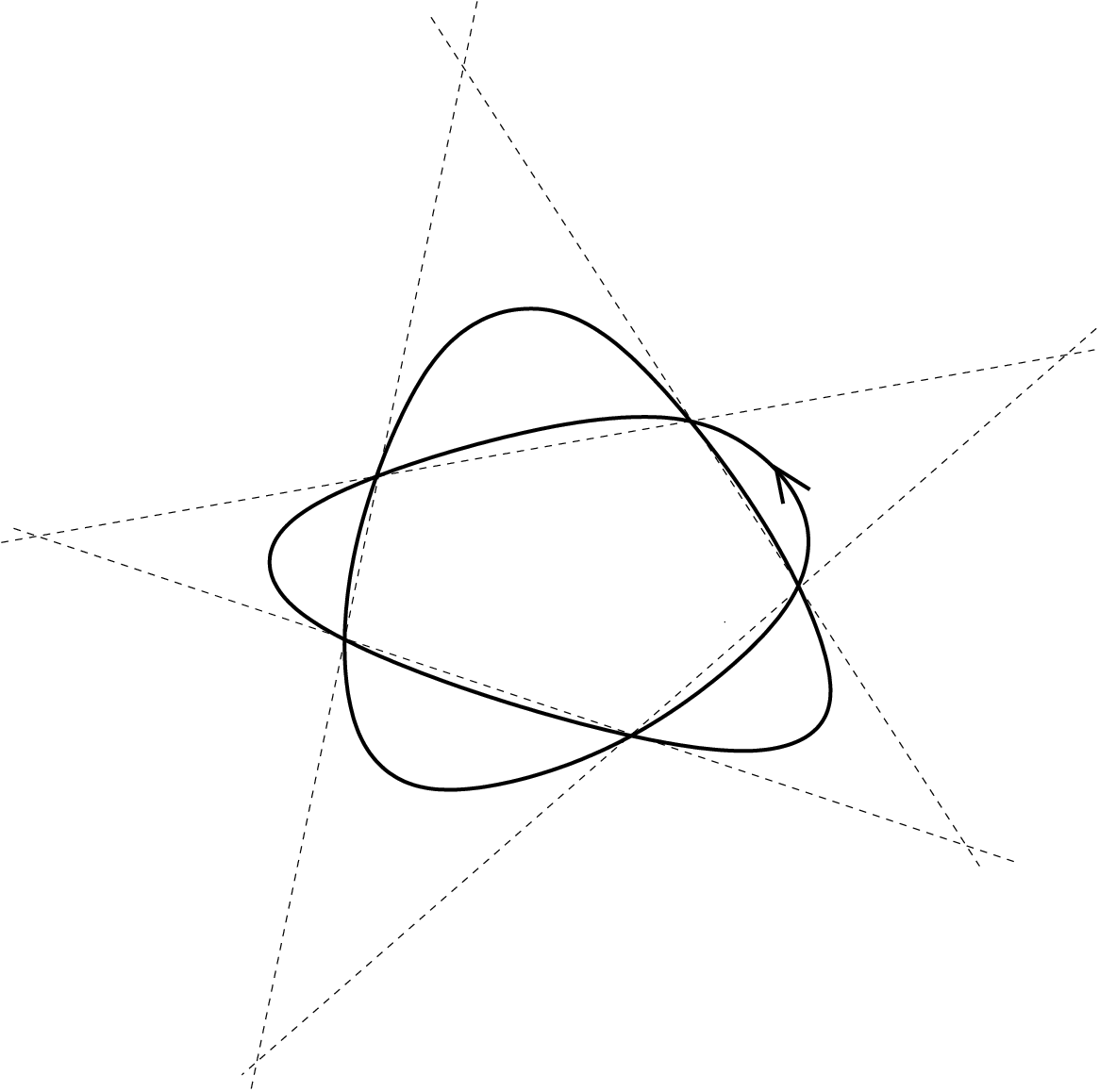}
\end{center}
\caption{A star is trapped in a union of triangles.}
\label{fig:starpoly}
\end{figure}

The only possible star from which a Reidemeister move of type III
is possible is the one shown in figure \ref{fig:reidstar} ($k = 1$):
indeed, a Reidemeister move of type III is quite impossible if the
curve does not form a combinatorial triangle.
In order to see that the only possible Reidemeister moves of type II
are those indicated in figure \ref{fig:reidstar}, notice that
if $\gamma$ is a star, its image is trapped in the union of triangles
shown in figure \ref{fig:starpoly} (where straight lines indicate
geodesics in the sphere).

For $\gamma \in \cL_{+1}^{(1)} \cap \cT_0$,
let $K(\gamma) \subset \Ss^2$ be the closure of the region
positively surrounded by $\gamma$.
The set $K(\gamma)$ is shaded in Figure \ref{fig:star}.
Clearly, $K(\gamma)$ is a convex set.
The definition of $K(\gamma)$ can be continuously extended to $\cT_0$:
for $\gamma \in \cT_1$, $K(\gamma)$ consists of the triple point only.
Let $k(\gamma) \in K(\gamma)$ be the baricenter of $K(\gamma)$
(recall that in order to find the baricenter of a convex
subset of $\Ss^2$ we first find its baricenter in $\RR^3$
and then radially project it onto the sphere).
For $\gamma \in \cT_1$, $k(\gamma)$ is the triple point.
Figure \ref{fig:starpoly} shows that $\gamma(t) \ne \pm k(\gamma)$
for any $\gamma \in \interior(\cT_0)$ and any $t$.

We construct a homotopy $H: [0,1] \times \interior(\cT_0)
\to \interior(\cT_0)$,
$H(0,\gamma) = \gamma$, $H(1,\cdot)$ constant equal to $\nu_2$.
Given $\gamma$, let $v_1, v_2, v_3$ be the only positively oriented
orthonormal basis with $v_3 = k(\gamma)$, $\gamma(0)$ in the plane
spanned by $v_1$ and $v_3$. Reparametrizing, we may write
\[ \gamma(t) = \frac{1}{\sqrt{1 + (u(t))^2}}\;
\left( (\cos(4\pi t)) v_1 + (\sin(4\pi t)) v_2 + u(t) v_3 \right). \]
The condition for such a curve to be locally convex
is that $u''(t) + 16\pi^2 u(t) > 0$.
Set $u(s,t) = s u(t) + (1-s)$ and (up to base point)
\[ H(s,\gamma)(t) = \frac{1}{\sqrt{1 + (u(s,t))^2}}\;
\left( (\cos(4\pi t)) v_1 + (\sin(4\pi t)) v_2 + u(s,t) v_3 \right). \]
By the linearity of the above condition, all such curves are locally convex.

We construct the universal cover of $\cT_1$.
Let $\tilde\cT_1$ be the set of pairs $(\tilde\gamma,\tilde t_0)$
where $\tilde\gamma: \RR \to \Ss^2$ is
a $1$-periodic locally convex function
with $\fF_{\tilde\gamma}(0) = I$,
the restriction $\gamma = \tilde\gamma|_{[0,1]}$ belongs to $\cT_1$
and $\tilde t_0 \in \RR$ is a triple point, i.e.,
there exist $\tilde t_1, \tilde t_2 \in RR$,
$\tilde t_0 < \tilde t_1 < \tilde t_2 < \tilde t_0 + 1$,
$\tilde\gamma(\tilde t_0) = \tilde\gamma(\tilde t_1) =
\tilde\gamma(\tilde t_2)$.
The projection $\Pi: \tilde\cT_1 \to \cT_1$ takes
$(\tilde\gamma,\tilde t_0) \in \tilde\cT_1$ to 
$\gamma = \tilde\gamma|_{[0,1]} \in \cT_1$.
This is a covering map by construction;
the group of deck transformations is isomorphic to $\ZZ$,
spanned by
$(\tilde\gamma,\tilde t_0) \mapsto (\tilde\gamma,\tilde t_1)$
(where $\tilde t_1$ is defined as above).

We claim that $\tilde\cT_1$ is contractible.
A homotopy $H: [0,1] \times \tilde\cT_1 \to \tilde\cT_1$ 
taking $\tilde\cT_1$ to a point starts with, for $s \in [0,1/4]$,
\[ H(s,\tilde\gamma,\tilde t_0) =
(\tilde\gamma_s,(1 - 4s)\tilde t_0), \quad
\tilde\gamma_s(t) =
(\fF_{\tilde\gamma(4s\tilde t_0)})^{-1} \tilde\gamma(t+4s\tilde t_0). \]
This defines a deformation retract from $\tilde\cT_1$ to $\cF_2$,
the set of flowers with $3$ petals (\cite{S1}).
We now use the interval $s \in [1/4,1/2]$
to reparametrize our curves so that for $s = 1/2$ the triple
point of $H(s,\tilde\gamma,\tilde t_0)$ will be
$t_0 = 0$, $t_1 = 1/3$, $t_2 = 2/3$.
Next, for $s \in [1/2,3/4]$,
set $H(s,\tilde\gamma,\tilde t_0) = (H(1/2,\tilde\gamma,\tilde t_0))^U$,
$U = U(\tilde\gamma,\tilde t_0) \in \cUpu$, so that
\[ \fF_{H(3/4,\tilde\gamma,\tilde t_0)}(1/3) = \begin{pmatrix} 1 & 0 & 0 \\
0 & -\frac{1}{2} & \frac{\sqrt{3}}{2} \\
0 & -\frac{\sqrt{3}}{2} & -\frac{1}{2} \end{pmatrix}, \quad
\fF_{H(3/4,\tilde\gamma,\tilde t_0)}(2/3) = \begin{pmatrix} 1 & 0 & 0 \\
0 & -\frac{1}{2} & -\frac{\sqrt{3}}{2} \\
0 & \frac{\sqrt{3}}{2} & -\frac{1}{2} \end{pmatrix}. \]
Finally, the loops $H(3/4,\tilde\gamma,\tilde t_0)|_I$,
$I = [0,1/3], [1/3,2/3], [2/3,1]$, are convex:
we use $s \in [3/4,1]$ to deform them to some fixed loop.
This completes the proof of the claim.
Thus $\pi_1(\cT_1) = \ZZ$ and $\cT_1$ has a contractible universal cover,
proving that $\cT_1$ is homotopically equivalent to $\Ss^1$
and completing the proof of the lemma.
\qed

\section{Eggs}
\label{sect:egg}

\begin{figure}[ht]
\begin{center}
\psfrag{(a)}{$(a)$}
\psfrag{(b)}{$(b)$}
\epsfig{height=35mm,file=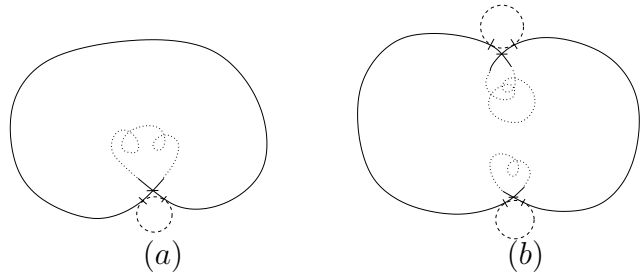}
\end{center}
\caption{Odd and even eggs.}
\label{fig:egg}
\end{figure}

A curve $\gamma \in \cL_I$ is an \textit{odd egg}
if there exists a reverse loop $(t_0,t_1)$
such that the image of $\gamma$ is contained in the closed
disk positively surrounded by the loop (Figure \ref{fig:egg}, (a)).
Notice that a star with a single double point is an odd egg.
A curve $\gamma \in \cL_I$ is an \textit{even egg}
if there exist two transversal double points
$(t_0,t_1)$ and $(t_2,t_3)$ with
$t_0 < t_1 < t_2 < t_3 < t_0 + 1$,
$\det(\gamma(t_0),\gamma'(t_0),\gamma'(t_1)) < 0$ and
$\det(\gamma(t_2),\gamma'(t_2),\gamma'(t_3)) < 0$
such that $\gamma|_{[t_1,t_2) \cup [t_3,t_0+1)}$ is injective
and the image of $\gamma$ is contained in the closed disk
positively surrounded by the above restriction (Figure \ref{fig:egg}, (b)).
Let $\cE \subset \cL_I$ be the set of all eggs (even and odd).

\begin{lemma}
\label{lemma:noeggs}
Let $K$ be a compact manifold and
$f: K \to \cL_I \smallsetminus (\cL_{-1,c} \cup \cT_0)$ a continuous map.
Then $f$ is homotopic in $\cL_I \smallsetminus (\cL_{-1,c} \cup \cT_0)$
to some map $\hat f: K \to \cL_I \smallsetminus
(\cL_{-1,c} \cup \cT_0 \cup \cE)$.
\end{lemma}


{\nobf Proof: }
Intuitively, we pull the creature out of the egg.

There exist open sets $B_1 \subset B_2 \subset K$,
$\overline{B_1} \subset B_2$, 
such that:
\begin{enumerate}[(a)]
\item{if $p \notin B_1$ then $f(p)$ is not an odd egg;}
\item{there exist functions $t_0, t_1: B_2 \to \Ss^1$ 
with $(f(p),t_0(p),t_1(p)) \in \cW$ for all $p \in B_2$;}
\item{if $f(p)$ is an odd egg then $f(p)(t_0(p))$ and $f(p)(t_1(p))$
are approximately equal to the extrema of the reverse loop
in the definition of odd eggs;}
\end{enumerate}
The functions $t_0$ and $t_1$ are indicated in Figure \ref{fig:egg}(a).
Let $\phi: K \to [0,1]$ be a continuous function with
$\phi|_{B_1} = 1$, $\phi|_{K \smallsetminus B_2} = 0$.
Define 
\[ H(s,p) = \begin{cases}
f(p), & p \notin B_2,\\
\Delta^\sharp(s \phi(p), f(p), t_0(p), t_1(p)), & p \in B_2.
\end{cases} \]
This opens all odd eggs.

In order to get rid of the even eggs the construction is similar
but with a harmless subtlety.
We can easily define $B_3 \subset B_4 \subset K$,
$\overline{B_3} \subset B_4$, 
such that if $p \notin B_3$ then $f(p)$ is not an even egg.
In $B_4$ is not simply connected,
it is not clear, however, that continuous functions
$t_0, t_1, t_2, t_3: B_4 \to \Ss^1$ can be defined
since the two double points in the shell of the egg can trade places.
Define therefore 
$t_0, t_1, t_2, t_3: \tilde B_4 \to \Ss^1$,
where $\tilde B_4$ is an appropriate double cover of $B_4$.
Finally, define
\[ H(s,p) = \begin{cases}
f(p), & p \notin B_4,\\
\Delta^\sharp(s \phi(p),
\Delta^\sharp(s \phi(p), f(p), t_0(\tilde p), t_1(\tilde p)),
t_2(\tilde p), t_3(\tilde p)), & p \in B_4,
\end{cases} \]
where $\tilde p \in \tilde B_4$ is one of the two lifts of $p$;
Lemma \ref{lemma:disjDeltasharp} guarantees that both choices of 
$\tilde p$ obtain the same value for $H$.
\qed

\begin{lemma}
\label{lemma:gammahasloop}
If $\gamma \in \cL_I^{(3)} \smallsetminus (\cL_{-1,c} \cup \cT_0 \cup \cE)$
then $\gamma$ has at least one direct loop;
if $\gamma \in \cL_I^{(2)} \smallsetminus (\cL_{-1,c} \cup \cT_0 \cup \cE)$
then $\gamma$ has at least two direct loops.
\end{lemma}

\begin{figure}[ht]
\begin{center}
\psfrag{tast}{$t_\ast$}
\psfrag{ta}{$t_a$}
\psfrag{tb}{$t_b$}
\psfrag{tta}{$\tilde t_a$}
\psfrag{ttb}{$\tilde t_b$}
\epsfig{height=25mm,file=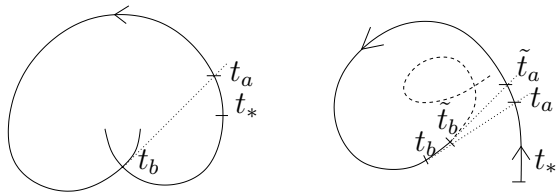}
\end{center}
\caption{A configuration which obtains a direct loop.}
\label{fig:nautilus}
\end{figure}

{\nobf Proof: }
We first consider the case $\gamma \in \cL_I^{(1)} \smallsetminus \cL_{-1,c}$.
Take $t_\ast \in \Ss^1$. Let
\[ t_1 = \sup \{ t \in \Ss^1 \;|\;
\gamma|_{[t_\ast,t]} \textrm{ is injective} \}. \]
There exists a unique $t_0 \in [t_\ast,t_1)$ with
$\gamma(t_0) = \gamma(t_1)$.
The desired loop is $(t_0,t_1)$.
Still in $\cL_I^{(1)} \smallsetminus \cL_{-1,c}$,
we prove the existence of a direct loop.
Let $(t_c,t_b)$ be a reverse loop.
As in Figure \ref{fig:nautilus} (a),
draw a geodesic tangent to $\gamma$ at $\gamma(t_b)$:
the geodesic transversally intersects the image of $\gamma$
at $\gamma(t_a)$, $t_a \in (t_c,t_b)$.
Take $t_\ast \in (t_c,t_a)$:
we claim that the construction above obtains a direct loop.
More generally, assume the restriction of $\gamma$
to $[t_\ast,t_b]$ is as in Figure \ref{fig:nautilus} (b):
an injective function such that the geodesic tangent to
the image of $\gamma$ at $\gamma(t_b)$ meets
the image of $\gamma$ transversally at $\gamma(t_a)$, $t_a \in [t_\ast,t_b)$;
also, the image under $\gamma$ of $[t_a,t_b]$ plus
the segment of geodesic between $\gamma(t_b)$ and $\gamma(t_a)$
form the boundary of a convex closed disk $D(t_b) \subset \Ss^2$.
Then, as $\tilde t_b$ increases starting from $t_b$ the above condition
in preserved (and $D(\tilde t_b)$ becomes smaller)
until $\gamma(\tilde t_a) = \gamma(\tilde t_b)$,
obtaining a direct loop.

The same construction and argument may be applied with time reversed:
it follows that the only curves in $\cL^{(1)}_I \smallsetminus \cL_{-1,c}$
with a unique direct loop are those for which both constructions
(original and with reversed time) lead to the same direct loop.
Thus, the only curves with a unique direct loop are odd eggs.

\begin{figure}[ht]
\begin{center}
\psfrag{gammat2}{$\gamma(t_2)$}
\epsfig{height=25mm,file=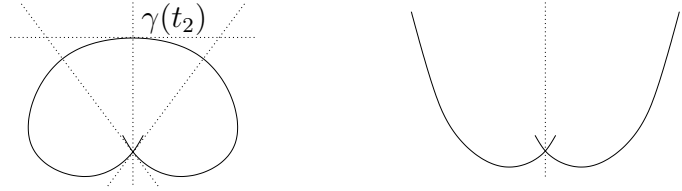}
\end{center}
\caption{A reverse loop and its image under a projective transformation.}
\label{fig:revloop}
\end{figure}

If $\gamma$ has no self-tangencies of odd order,
perturb it near each self-tangency
so as to destroy the self-tangency without creating new self-intersections.
Consider a loop of the modified curve $\hat\gamma$:
we claim that the same double point is a loop for the original curve $\gamma$.
It suffices to show that self-tangencies can not be created within loops.
For direct loops this follows from the convexity of the restriction.
For a reverse loop $(t_0,t_1)$,
take the two geodesics tangent to $\gamma$ at $t_0$ and $t_1$
as in Figure \ref{fig:revloop} and a point $\gamma(t_2)$
between the intersections of these geodesics with the image of $\gamma$.
Consider a projective transformation taking the tangent geodesic to $\gamma$
at $\gamma(t_2)$ to infinity and the geodesic joining
$\gamma(t_0) = \gamma(t_1)$ with $\gamma(t_2)$.
The image of $\gamma$ under this projective transformation
is the graph of a function, completing the proof in this case.

If $\gamma$ has a single self-tangency $(t_0,t_1)$ of order $3$
start by perturbing a neighborhood of the self-tangency
in such a way as to create exactly one transversal double point.
Any direct loop of $\hat\gamma$ except for $(t_0,t_1)$ (if it is simple)
yields a direct loop of $\gamma$.
Thus, the only situation where $\gamma$ does not have a direct loop
is if $(t_0,t_1)$ is the unique direct loop of $\hat\gamma$;
this completes the proof.
\qed 

Two loops $(t_0,t_1)$ and $(\tilde t_0, \tilde t_1)$ are \textit{disjoint}
if the intervals $[t_0,t_1] \subset \Ss^1$
and $[\tilde t_0, \tilde t_1] \subset \Ss^1$ are disjoint.
Notice that this does not mean that the images of the intervals
under $\gamma$ are disjoint.

\begin{figure}[ht]
\begin{center}
\epsfig{height=25mm,file=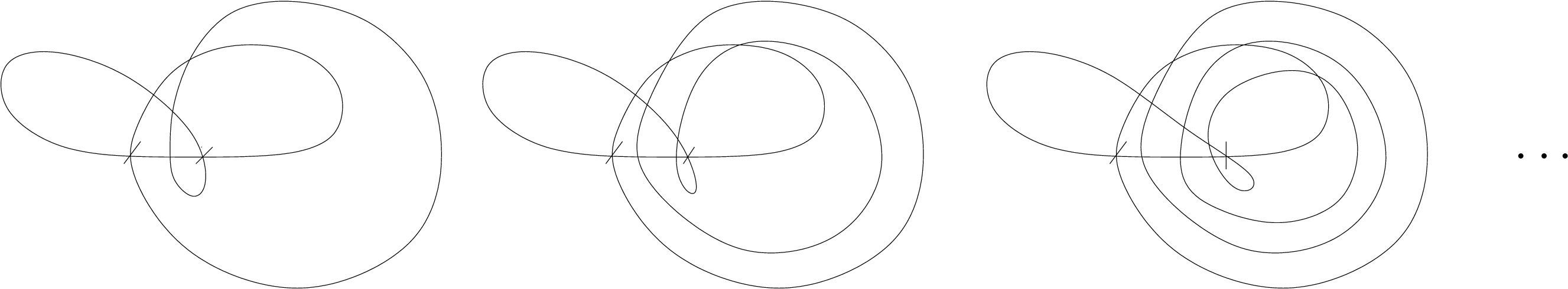}
\end{center}
\caption{Pseudo-eggs}
\label{fig:pseudoegg}
\end{figure}

A curve $\gamma \in \cL^{(1)}_I \smallsetminus
(\cL_{-1,c} \cup \cT_0 \cup \cE)$
is a \textit{pseudo-egg} if it belongs to the connected component
of one of the curves in the infinite family
indicated in Figure \ref{fig:pseudoegg}.
More precisely, $\gamma$ admits two non disjoint direct loops
$(t_0,t_1)$ and $(t_2,t_3)$, $t_0 < t_2 < t_1 < t_3 < t_0 + 1$,
and the restriction of $\gamma$ to $[t_3,t_0+1]$ is injective.
Let $\cEp$ be the set of eggs and pseudo-eggs.

\begin{lemma}
\label{lemma:nopseudoeggs}
Let $K$ be a compact manifold and
$f: K \to \cL_I \smallsetminus
(\cL_{-1,c} \cup \cT_0 \cup \cE)$ a continuous map.
Then $f$ is homotopic in $\cL_I \smallsetminus
(\cL_{-1,c} \cup \cT_0 \cup \cE)$ to some map
$\hat f: K \to \cL_I \smallsetminus
(\cL_{-1,c} \cup \cT_0 \cup \cEp)$.
\end{lemma}

\begin{figure}[ht]
\begin{center}
\psfrag{or}{or}
\epsfig{height=25mm,file=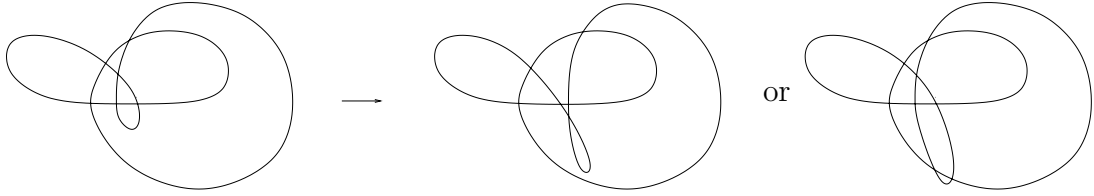}
\end{center}
\caption{How to get rid of pseudo-eggs}
\label{fig:nomopseudoegg}
\end{figure}

{\nobf Proof: }
Consider the innermost loop in the spiral and pull it
(using $\Delta^\sharp$) as in Figure \ref{fig:nomopseudoegg}:
this will either do a Reidemeister move of type III
or of type II, in either case destroying the pseudo-egg.
\qed

\begin{lemma}
\label{lemma:gammahasdisjloops}
Consider a curve $\gamma \in \cL^{(1)}_I \smallsetminus
(\cL_{-1,c} \cup \cT_0 \cup \cEp)$.
For every direct loop $(t_0,t_1)$ of $\gamma$
there is a disjoint direct loop $(t_2,t_3)$.
For every reverse loop $(t_0,t_1)$ of $\gamma$
there is a disjoint (reverse or direct) loop $(t_2,t_3)$.
Furthermore, given two non-disjoint direct loops $(t_0,t_1)$ and $(t_2,t_3)$,
$\gamma$ admits a loop $(t_4,t_5)$
which is disjoint from both $(t_0,t_1)$ and $(t_2,t_3)$.
\end{lemma}

\begin{figure}[ht]
\begin{center}
\psfrag{(a)}{(a)}
\psfrag{(b)}{(b)}
\epsfig{height=30mm,file=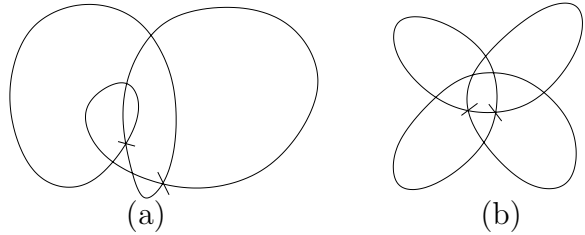}
\end{center}
\caption{Non-disjoint loops $(t_0,t_1)$ and $(t_2,t_3)$.}
\label{fig:counterhasdisjloops}
\end{figure}

The second part does not always hold
if $(t_0,t_1)$ and $(t_2,t_3)$ are reverse:
see Figure \ref{fig:counterhasdisjloops}, (a).
In the last claim, the new loop in case (a) may be reverse:
see Figure \ref{fig:counterhasdisjloops}, (b).

{\nobf Proof: }
Assume $\gamma \in \cL^{(1)}_I$ so that every self-intersection is transversal.
Let $(t_0,t_1)$ be a loop.
It $(t_1,t_0+1)$ is a loop then $\gamma$ is a star.
Otherwise, take $t_\ast = t_1 + \epsilon$
($\epsilon > 0$, $\epsilon$ small).
As in Lemma \ref{lemma:gammahasloop},
let $t_3$ be the smallest $t > t_\ast$ such that $\gamma|_{[t_\ast,t]}$
is not injective and let $t_2 \in [t_\ast,t_3)$ be such that
$\gamma(t_2) = \gamma(t_3)$, so that $(t_2,t_3)$ is a loop.
Since $(t_1,t_0+1)$ is not a loop, $t_3 < t_0 + 1$ and
$(t_2,t_3)$ is disjoint from $(t_0,t_1)$, as required.

\begin{figure}[ht]
\begin{center}
\psfrag{gt0gt1}{$(t_0,t_1)$}
\psfrag{gt2gt3}{$(t_2,t_3)$}
\epsfig{height=25mm,file=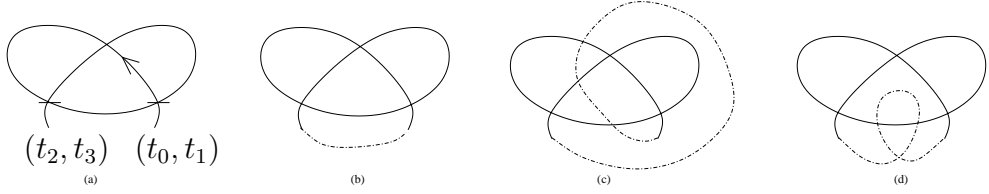}
\end{center}
\caption{Two non-disjoint loops.}
\label{fig:nondisjointloops}
\end{figure}

Let $(t_0,t_1)$, $(t_2,t_3)$ be two non-disjoint loops.
The case $t_0 = t_2$, $t_1 = t_3$ was discussed in the previous paragraph;
$t_2 = t_1$, $t_3 = t_0 + 1$ implies that $\gamma$ is a star.
We may therefore assume $t_0 < t_2 < t_1 < t_3 < t_0+1$,
as in Figure \ref{fig:nondisjointloops}.

Notice that $(\fF_\gamma(t_0))^{-1} \fF_\gamma(t_3) \in \cJ_2$.
If the restriction $\gamma|_{[t_3,t_0+1]}$ is convex then
$\gamma$ is a star
(Figure \ref{fig:nondisjointloops} (b)).
If the restriction is injective but not convex
then $\gamma$ is either an egg or a pseudo-egg
(Figure \ref{fig:nondisjointloops} (c)).
Finally, if the above restriction is not injective then a new loop
disjoint from $(t_0,t_1)$ and $(t_2,t_3)$ exists
(Figure \ref{fig:nondisjointloops} (d)).
\qed




\section{Proof of the main results}

\begin{lemma}
\label{lemma:pi1minusT}
The connected components of $\cL_I \smallsetminus (\cL_{-1,c} \cup \cT_0)$
are simply connected.
\end{lemma} 

{\nobf Proof: }
Let $f: \Ss^1 \to \cL_I \smallsetminus (\cL_{-1,c} \cup \cT_0)$
be a continuous map.
We prove that $f$ is homotopic to a point.
By Lemma \ref{lemma:noeggs} and transversality we may assume
$f: \Ss^1 \to \cL^{(2)}_I \smallsetminus (\cL_{-1,c} \cup \cT_0 \cup \cE)$
and that $f(p) \in \cL^{(1)}_I$ except for finitely many points
$p_1, \ldots, p_M \in \Ss^1$.
In case $M \le 2$ add a few points to the list to guarantee $M \ge 3$.
Reparametrize so that $p_k = \frac{k}{M} \in \Ss^1$.
By Lemma \ref{lemma:gammahasloop},
for every $p \in \Ss^1$ the curve $f(p)$ has a direct loop.
Let $V_k = (\frac{k-1}{M},\frac{k+1}{M})$ and use a direct loop of
$f(\frac{k}{M})$ to define $t_{0,k}, t_{1,k}: V_k \to \Ss^1$.
This is an open cover, but the loops are probably not disjoint.

Let $U_k = (\frac{3k-1}{3M},\frac{3k+1}{3M})$.
Let $p_\star = \frac{2k+1}{2M}$, $\gamma_\star = f(p_\star) \in \cL^{(1)}_I$,
$t_0 = t_{0,k}(p_\star)$, $t_1 = t_{1,k}(p_\star)$,
$t_2 = t_{0,k+1}(p_\star)$ and $t_3 = t_{1,k+1}(p_\star)$.
If the loops $(t_0,t_1)$ and $(t_2,t_3)$ of $\gamma_\star$ are equal
or disjoint, set $U_{k+1/2} = (\frac{k}{M},\frac{k+1}{M})$ and
$t_{\ast,k+1/2}(p) = t_{\ast,k}(p)$.
If these two loops are not disjoint, use Lemma \ref{lemma:gammahasdisjloops}
to obtain two loops $(t_4,t_5)$ and $(t_6,t_7)$ of $\gamma_\star$ such that:
\begin{enumerate}
\item{$(t_4,t_5)$ is disjoint from $(t_0,t_1)$;}
\item{$(t_6,t_7)$ is disjoint from $(t_2,t_3)$;}
\item{$(t_4,t_5)$ and $(t_6,t_7)$ are either equal or disjoint.}
\end{enumerate}
Set $U_{k+1/3} = (\frac{k}{M},\frac{3k+2}{3M})$ and
$U_{k+2/3} = (\frac{3k+1}{3M},\frac{k+1}{M})$;
use $(t_4,t_5)$ and $(t_6,t_7)$ to define
$(t_{0,k+1/3},t_{1,k+1/3})$
and $(t_{0,k+2/3},t_{1,k+2/3})$, respectively.
Thus, $f$ admits a cover by disjoint loops.
From Lemma \ref{lemma:disjcover}, $f$ is homotopic to $\nu_2 \ast f$.

On the other hand, each connected component of $\cI_I$
is simply connected and therefore $f$ is homotopic to a point in $\cI_I$.
By Proposition \ref{prop:S1} $f$ is homotopic to a point in $\cL_I$.
The homotopy can be constructed as follows:
consider an extension $\hat f: \DD^2 \to \cI_I$ of $f$;
use $\Delta^\sharp$ to pass from $f$ to $F_N \circ f$
and complete the homotopy with $F_N \circ \hat f$.
For sufficiently large $N$,
this homotopy remains in $\cL_I \smallsetminus (\cL_{-1,c} \cup \cT_0)$.
\qed

\begin{theo}
\label{theo:pi1}
Each connected component of $\cL_I$ is simply connected.
\end{theo}

{\nobf Proof: }
It is well known that $\cL_{-1,c}$ is contractible.
Since $\cT_0 \subset \cL_{+1}$,
we just proved that in Lemma \ref{lemma:pi1minusT}
that $\cL_{-1,n}$ is also simply connected.
Finally, we use Seifert-Van Kampen to compute $\pi_1(\cL_{+1})$.
Let $A = \cL_{+1} \smallsetminus \cT_0$ and
fatten $\cT_0$ a bit to obtain an open set $B$:
since $\cT_0$ is contractible and $\cT_1 = \partial\cT_0$
is a connected submanifold of codimension $1$,
$B$ is simply connected and $A \cap B$ is connected.
Thus, $\cL_{+1} = A \cup B$ is also simply connected.
\qed

\begin{lemma}
\label{lemma:filldisk}
Let $\gamma \in \cL^{(1)}_I \smallsetminus (\cL_{-1,c} \cup \cT_0 \cup \cEp)$.
Let $f: \DD^2 \to \cL_I \smallsetminus \cT_0$ be the function
constant equal to $\gamma$.
Consider a disjoint cover by loops $\cC$ of $f|_{\Ss^1}$
without two consecutive reverse loops.
Then there exists $\tilde f: \DD^2 \to \cL_I \smallsetminus \cT_0$
homotopic to $f$ with fixed boundary
and an extension $\tilde\cC$ of $\cC$
to a disjoint cover by loops of $\tilde f$.
\end{lemma}

\begin{figure}[ht]
\begin{center}
\psfrag{V0}{$V_0$}
\psfrag{Vi}{$V_i$}
\psfrag{tV}{$\tilde V$}
\psfrag{(a)}{(a)}
\psfrag{(b)}{(b)}
\psfrag{(c)}{(c)}
\epsfig{height=45mm,file=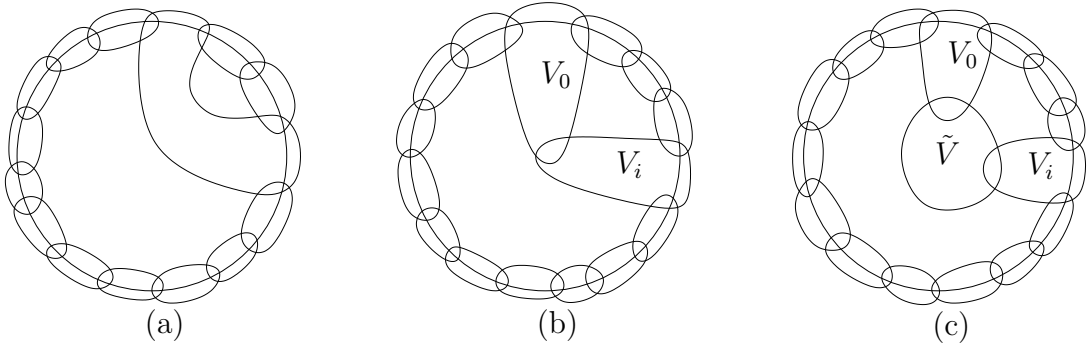}
\end{center}
\caption{Long cycles can be decomposed into short ones.}
\label{fig:nolongcycle}
\end{figure}

{\nobf Proof: }
Without loss of generality the open sets in the cover are intervals.
If two neighboring intervals have identical loops, fuse them;
may therefore assume without loss of generality that $\Ss^1$
is covered by a cycle of loops, where two adjacent loops are disjoint
and no two adjacent loops are both reverse.
More: if two non-adjacent intervals use the same loop,
the corresponding open set can be enlarged to cross the disk
(Figure \ref{fig:nolongcycle},(a)).
We may therefore assume all loops to be distinct.

If the cycle $7$ or more loops, it can be decomposed into cycles
of size $6$ or less.
Indeed, consider a cycle of length $n > 6$;
let $\ell_0$ be a direct loop in the cycle
and let $\ell_1, \ell_2, \ell_3, \ell_4$
be its neighbors (numbered clockwise).
Let $i = 3$ or $4$ such that $\ell_i$ is direct.
If $\ell_0$ and $\ell_i$ are disjoint, extend the open sets
$V_0$ and $V_i$ to intersect in the center,
thus subdividing the original cycle into one of length $i+1$
and another of length $n-i+1$ (Figure \ref{fig:nolongcycle},(b)).
If $\ell_0$ and $\ell_i$ are not disjoint,
use Lemma \ref{lemma:gammahasdisjloops} to obtain
a loop $\tilde\ell$ disjoint from both:
introduce an open set $\tilde V$ with associated loop $\tilde\ell$
in the middle of the disk, intersecting $V_0$ and $V_i$ only
among the originally defined loops.
This subdivides the original cycle into two cycles of lengths $i+2$
and $n-i+2$ (Figure \ref{fig:nolongcycle},(c)).

Cycles of $2$ or $3$ loops have all loops disjoint and therefore
admit a disjoint cover as they are, no homotopy needed.
We are left with considering cycles of length $4$, $5$ or $6$.

From now on we argue by contradiction,
searching for the shortest counterexample.
If two non-adjacent loops are disjoint, not both reverse,
construction (b) in Figure \ref{fig:nolongcycle} can be applied,
decomposing our cycle into two shorter ones.
We may therefore assume that non-adjacent loops not both reverse
are not disjoint.

We first consider the case $n=6$. If two opposite loops
$\ell_i$ and $\ell_{i+3}$ are both direct then
the construction in Figure \ref{fig:nolongcycle} (b) or (c)
can be applied, decomposing our cycle.
We may therefore assume that $\ell_i$ is direct for $i$ even
and reverse for $i$ odd.
Consider the intervals $I_i = [t_{0,i},t_{1,i}] \subset \Ss^1$:
$I_0$, $I_2$ and $I_4$ are pairwise neither disjoint not nested
(one contained in the other).
There are, up to reparametrization and permutation,
only two possibilities:
$I_0 = [0,3/6]$, $I_2 = [1/6,4/6]$, $I_4 = [-1/6,2/6]$ or
$I_0 = [0,3/6]$, $I_2 = [2/6,5/6]$, $I_4 = [-2/6,1/6]$.
In either case, $I_3$ must be disjoint from both $I_2$ and $I_4$
and neither disjoint nor nested with $I_0$, a contradiction.

\begin{figure}[ht]
\begin{center}
\psfrag{l0}{$\ell_0$}
\psfrag{l2}{$\ell_2$}
\psfrag{l3}{$\ell_3$}
\psfrag{V0}{$V_0$}
\psfrag{V1}{$V_1$}
\psfrag{V2}{$V_2$}
\psfrag{V3}{$V_3$}
\psfrag{V4}{$V_4$}
\epsfig{height=45mm,file=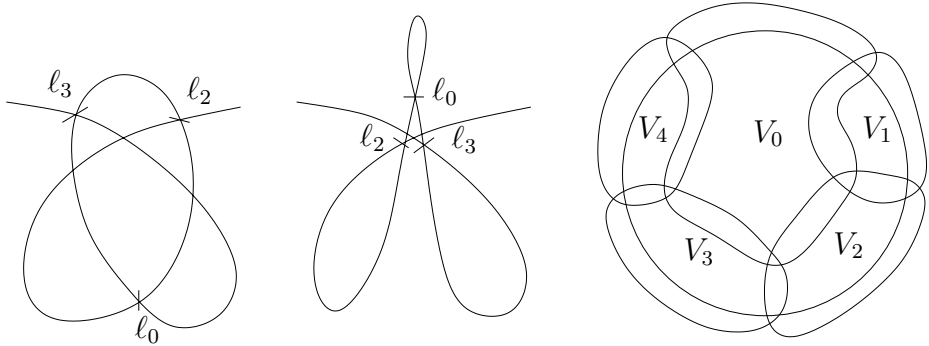}
\end{center}
\caption{The case $n=5$}
\label{fig:n5}
\end{figure}

Consider now the case $n=5$, with loops $\ell_0, \ldots, \ell_4$.
We may assume without loss of generality that
$\ell_0$, $\ell_2$ and $\ell_3$ are direct loops,
$\ell_2$ and $\ell_3$ disjoint but $\ell_0$ not disjoint from either.
We may again assume that
$I_0 = [1/10,4/10]$, $I_2 = [0,2/10]$, $I_3 = [3/10,5/10]$
and that the image of $[0,5/10]$ under $\gamma$ is
as in Figure \ref{fig:n5}.
Pull the loop $\ell_0$ (or, in other words,
apply $\Delta^\sharp(s,\gamma,1/10-\epsilon,4/10+\epsilon)$)
in the center of the disk $\DD^2$ to define $\tilde f$.
The loops $\ell_0, \ell_2, \ell_3$ survive in $\DD^2$
and near the center of $\DD^2$,
$\ell_0$ becomes disjoint from $\ell_2, \ell_3$.
The loops $\ell_1$ and $\ell_4$ were not affected
and remain disjoint from $\ell_0$.
We therefore have the disjoint cover in Figure \ref{fig:n5}.

\begin{figure}[ht]
\begin{center}
\psfrag{l0}{$\ell_0$}
\psfrag{l2}{$\ell_2$}
\psfrag{V0}{$V_0$}
\psfrag{V1}{$V_1$}
\psfrag{V2}{$V_2$}
\psfrag{V3}{$V_3$}
\epsfig{height=45mm,file=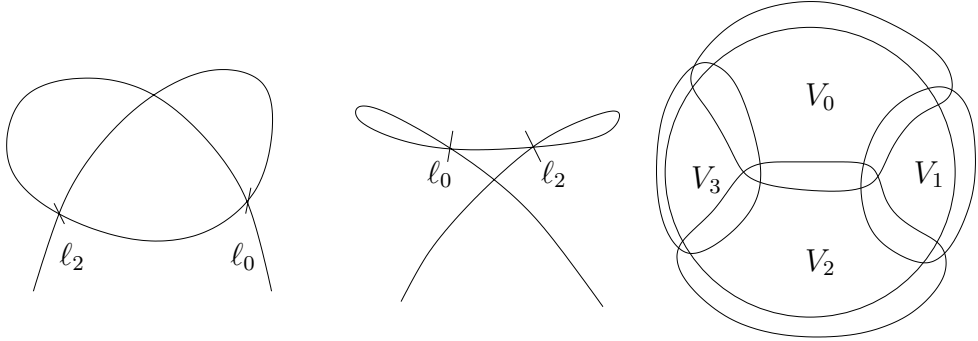}
\end{center}
\caption{The case $n=4$}
\label{fig:n4}
\end{figure}

Finally, for $n=4$, we may assume that $\ell_0$ and $\ell_2$ are direct.
As in Figure \ref{fig:n4}, pulling either $\ell_0$ or $\ell_2$
(or both) makes them disjoint.
The disjoint cover in Figure \ref{fig:n4} completes the proof.
\qed

\begin{lemma}
\label{lemma:pi2minusT}
Every continuous map $f: \Ss^2 \to \cL_I \smallsetminus
(\cL_{-1,c} \cup \cT_0)$
is homotopic to $\nu_2 \ast f$.
\end{lemma}

{\nobf Proof: }
It is enough to prove that $f$ is homotopic to some $\tilde f$
which admits a disjoint cover.
By Lemmas \ref{lemma:noeggs}, \ref{lemma:nopseudoeggs}
and transversality we may assume
$f: \Ss^1 \to \cL^{(3)}_I \smallsetminus (\cL_{-1,c} \cup \cT_0 \cup \cEp)$.
Again by transversality, we may assume that there exists
a finite triangulation of $\Ss^2$ such that
if $f(p) \notin \cL^{(2)}_I$ then $p$ is a vertex of the triangulation and
if $f(p) \notin \cL^{(1)}_I$ then $p$ belongs to an edge.
By Lemma \ref{lemma:gammahasloop}, there exist open neighborhoods $V_i$
of the vertices and direct loops $t_{0,i}, t_{1,i}: V_i \to \Ss^1$.
As in the proof of Lemma \ref{lemma:pi1minusT},
there exists a disjoint cover of each edge by loops
without two consecutive reverse loops.
It remain to fill in the faces:
this is precisely what Lemma \ref{lemma:filldisk} does.
\qed

\begin{theo}
\label{theo:pi2}
We have $\pi_2(\cL_{+1}) = \ZZ^2$, $\pi_2(\cL_{-1,n}) = \ZZ$,
$H^2(\cL_{+1};\ZZ) = \ZZ^2$ and $H^2(\cL_{-1,n};\ZZ) = \ZZ$.
\end{theo}

{\nobf Proof: }
By the previous lemma, $f: \Ss^2 \to \cL_{-1,n}$ is homotopic to a point
in $\cL_{-1,n}$ if and only if it is homotopic to a point in $\cI_{-1,n}$.
In other words, the inclusion $\cL_{-1,n} \subset \cI_{-1}$
induces an isomorphism between $\pi_2(\cL_{-1,n})$ and $\pi_2(\cI_{-1}) = \ZZ$.
By Hurewicz theorem, inclusion also yields an isomorphism
between $H^2(\cL_{-1,n};\ZZ)$ and $H^2(\cI_{-1};\ZZ) = \ZZ$.
In other words, $H^2(\cL_{-1,n};\ZZ)$ is generated by $\xx$.

Similarly, $H^2(\cL_{+1} \smallsetminus \cT_0;\ZZ)$ is generated by $\xx$.
Use the normal bundle to $\cT_1$ to define open sets $A$ and $B$,
$\cT_0 \subset A$, $\cL_{+1} \smallsetminus \cT_0 \subset B$
such that the above inclusions and $\cT_1 \subset A \cap B$
are homotopy equivalences.
Write the Mayer-Vietoris sequence (coefficients in $\ZZ$):
\[  H^1(A) \oplus H^1(B) \to H^1(A \cap B) \to H^2(A \cup B) \to
H^2(A) \oplus H^2(B) \to H^2(A \cap B)  \]
We know that $H^1(A) = H^1(B) = 0$, $H^1(A \cap B) = \ZZ$,
$H^2(A) = 0$, $H^2(B) = \ZZ$, $H^2(A \cap B) = 0$.
Thus $H^2(A \cup B) = \ZZ^2$.
In the proof of Lemma \ref{lemma:trefoil}
we saw a geometric description of the universal cover
and therefore of the generator of $H^1(A \cap B)$:
appying the map $H^1(A \cap B) \to H^2(A \cup B)$
to this generator obtains $\ff_2 \in H^2(\cL_{+1})$,
the intersection number with $\cF_2 \subset \cT_1$.
Thus, $H^2(\cL_{+1}$ is generated by $\xx$ and $\ff_2$.
Again by Hurewicz theorem, $\pi_2(\cL_{+1}) = \ZZ^2$
is generated by $\bfg_2$ and $\nu_2 \ast \bfg_2$.
\qed

\section{Final remarks}

It may be possible to carry further methods used in this paper
to compute $H^\ast(\cL_I)$, but the kind of case-by-case analysis
in Section \ref{sect:egg} would have to be replaced by something
less accidental.
We hope to do this in \cite{Saldanha3} to prove
that the classes $\xx^n$ and $\ff_{2n}$
are generators of $H^\ast(\cL_{\pm 1})$
and that $\cL_{+1}$ and $\cL_{-1,n}$ have the homotopy type of
$\Omega\Ss^3 \vee \Ss^2 \vee \Ss^6 \vee \Ss^{10} \vee \cdots$ and
$\Omega\Ss^3 \vee \Ss^4 \vee \Ss^8 \vee \Ss^{12} \vee \cdots$,
respectively.


\bigskip\bigskip\bigbreak

{

\parindent=0pt
\parskip=0pt
\obeylines

Nicolau C. Saldanha, PUC-Rio
saldanha@puc-rio.br; http://www.mat.puc-rio.br/$\sim$nicolau/



\smallskip

Departamento de Matem\'atica, PUC-Rio
R. Marqu\^es de S. Vicente 225, Rio de Janeiro, RJ 22453-900, Brazil

}


\bigskip

\begin{thebibliography}{[10]}

\bibitem{BZ}{ Burde, G. and Zieschang, H.,
{\sl Knots}, Walter de Gruyter, Berlin, 1985.}
\bibitem{Little}{ Little, J. A.,
{\sl Nondegenerate homotopies of curves on the unit 2-sphere},
J. Differential Geometry, 4, 339-348, 1970.}
\bibitem{Saldanha3}{Saldanha, N.,
{\sl The homotopy type of spaces of locally convex curves in the sphere},
in preparation.}
\bibitem{S1}{Saldanha, N.,
{\sl The cohomology of spaces of locally convex curves in the sphere --- I},
preprint, arXiv:0905.2111v1.}
\bibitem{SaSha}{Saldanha, N. and Shapiro, B.,
{\sl Spaces of locally convex curves in $S^n$ and
combinatorics of the group $B^+_{n+1}$},
Journal of Singularities, 4, 1-22, 2012.}
\bibitem{SK}{Shapiro, B. and Khesin, B.,
{\sl Homotopy classification of nondegenerate quasiperiodic curves
on the $2$-sphere},
Publ. Inst. Math. (Beograd) 66(80), 127-156, 1999.}
\bibitem{ShapiroM}{Shapiro, M.,
{\sl Topology of the space of nondegenerate curves},
Math. USSR, 57, 106-126, 1993.}

\end{thebibliography}
\end{document}